\documentclass[11pt]{article}
\usepackage{graphicx}
\usepackage{textcomp}
\usepackage{amsmath,amsfonts,amscd,amsthm,amssymb}
\usepackage{epsfig}
\usepackage{makeidx}
\usepackage{color} 

\makeindex
\newcommand{\ncom}{\newcommand}
\newcommand{\be}{\begin{eqnarray*}}
\newcommand{\ee}{\end{eqnarray*}}
\newcommand{\ben}{\begin{eqnarray}}
\newcommand{\een}{\end{eqnarray}}

\newcommand{\tl}{\tilde}

\newcommand{\pr}{\partial}
\newcommand{\nab}{\nabla}
\newcommand{\bu}{{\bf u}}
\newcommand{\bU}{{\bf U}}

\newcommand{\bw}{{\bf w}}
\newcommand{\e}{{\bf e}}
\newcommand{\bs}{{\bf \sigma}}

\newcommand{\bJ}{{\bf J}}

\newcommand{\bv}{{\bf v}}

\newcommand{\bzeta}{\mbox{\boldmath $\zeta$}}

\newcommand{\bphi}{\mbox{\boldmath $\phi$}}
\newcommand{\brho}{\mbox{\boldmath $\rho$}}

\newcommand{\bxi}{\mbox{\boldmath $\xi$}}
\newcommand{\bH}{{\bf {H}}}

\newcommand{\bL}{\bf {L}}

\newcommand{\f}{{\bf {f}}}

\newtheorem{tdf}{Theorem}[section]

\newtheorem{ldf}{Lemma}[section]
 
\ncom{\ul}{\underline}
\ncom{\beq}{\begin{equation}}
\ncom{\eeq}{\end{equation}}
\ncom{\bea}{\begin{eqnarray*}}
\ncom{\eea}{\end{eqnarray*}}
\ncom{\beqa}{\begin{eqnarray}}
\ncom{\eeqa}{\end{eqnarray}}
\ncom{\nno}{\nonumber}
\ncom{\non}{\nonumber}
\ncom{\ds}{\displaystyle}
\ncom{\half}{\frac{1}{2}}
\ncom{\mbx}{\makebox{.25cm}}
\ncom{\hs}{\mbox{\hspace{.25cm}}}
\ncom{\rar}{\rightarrow}
\ncom{\Rar}{\Rightarrow}
\ncom{\noin}{\noindent}
\ncom{\sz}{\scriptsize}
\ncom{\Sgm}{\Sigma}
\ncom{\psgm}{\sigma^{\prime}}
\ncom{\dt}{\delta}
\ncom{\Dt}{\Delta}
\ncom{\lmd}{\lambda}
\ncom{\Lmd}{\Lambda}
\ncom{\Th}{\Theta}
\ncom{\eps}{\epsilon}
\ncom{\pcc}{\stackrel{P}{>}}
\ncom{\lp}{\stackrel{L_{p}}{>}}
\ncom{\sspan}{{\rm\,span}}
\ncom{\re}{{\rm Re\,}}
\ncom{\im}{{\rm Im\,}}
\ncom{\sgn}{{\rm sgn\,}}
\ncom{\ba}{\begin{array}}
\ncom{\ea}{\end{array}}
\ncom{\integ}[4]{\int_{#1}^{#2}\,{#3}\,d{#4}}
\ncom{\vspan}[1]{{{\rm\,span}\{ #1 \}}}
\ncom{\dm}[1]{ {\displaystyle{#1} } }
\ncom{\ri}[1]{{#1} \index{#1}}

\newtheorem{example}{Example}[section]

\newtheoremstyle
    {remarkstyle}
    {}
    {11pt}
    {}
    {}
    {\bfseries}
    {:}
    {     }
    {\thmname{#1} \thmnumber{#2} }

\theoremstyle{remarkstyle}

\begin{document}
\title
 {On a three level two-grid finite element method for the 2D-transient Navier-Stokes equations}
\author
 {Saumya Bajpai\\
TIFR Centre for Applicable Mathematics, \\
Post Bag No. 6503, GKVK Post Office,
Sharada Nagar, \\
Chikkabommasandra,
Bangalore 560065, India,\\
 Amiya K. Pani \\
 Department of Mathematics, Industrial Mathematics Group, \\
Indian Institute of Technology Bombay, Powai, Mumbai-400076, India
}
 \maketitle 
\abstract{In this paper, an error analysis of a three steps two level Galekin finite 
element method for the two dimensional transient Navier-Stokes equations is discussed.  First of all, the problem is 
discretized in spatial direction by employing finite element method on a coarse mesh $\mathcal{T} _H$ with mesh size $H$. 
Then, in step two, the nonlinear system is linearized around the coarse grid solution, say, $u_H$, which is similar 
to Newton's type iteration and the resulting linear system 
is solved on a finer mesh $\mathcal{T}_h$ with mesh size $h$. In step three, a correction is obtained 
through solving a linear problem on the finer mesh and an updated final solution is derived. Optimal error 
estimates  in $L^{\infty}({\bf L}^2)$-norm, when $h=\mathcal{O} (H^{2-\delta})$ and 
in $L^{\infty}(\bH^1)$-norm, when $h=\mathcal{O}(H^{4-\delta})$ for the velocity and in 
$L^{\infty}(L^2)$-norm,  when $h=\mathcal{O}(H^{4-\delta})$ for the pressure are established for 
arbitrarily small $\delta>0$. Further, under uniqueness assumption, 
these estimates are proved to be valid uniformly in time. Then based on backward Euler method, a 
completely discrete scheme is analyzed and {\it a priori} error estimates are derived. Finally, the paper 
is concluded with some numerical experiments.}


\section{Introduction}
\setcounter{equation}{0}
Consider the 2D-transient Navier-Stokes system:
\begin{eqnarray}
\label{1.1}
\frac {\partial \bu}{\partial t}-\nu \Delta \bu+ \bu\cdot \nabla \bu+\nabla p
={\bf f}({\bf x},t)\,\,\, {\bf x}\in \Omega ,\,\,t>0,
\end{eqnarray}
and incompressibility condition
\begin{eqnarray}
\label{1.2}
\nabla \cdot \bu=0\,\,\,{\bf x}\in \Omega,\,t>0,
\end{eqnarray}
with initial and boundary conditions
\begin{eqnarray}
\label{1.3}
\bu({\bf x},0)= \bu_0 \;\;\;\mbox {in}\;\Omega,\;\;\;\;\; \bu=0\;\; \;
\mbox {on}\; \partial \Omega,\; t\ge 0,
\end{eqnarray}
where, $\Omega$ is a bounded and convex polygonal domain in $\mathbb{R}^2$ with boundary $\partial \Omega$ and $f$ is 
given external force. Here, $\bu=\bu({\bf x},t)$ is the velocity vector, $p=p({\bf x},t)$ denotes the pressure 
and $\nu>0$ is the kinematic coefficient of viscosity.

In this article, a three level two-grid finite element Galerkin method for the problem (\ref{1.1})-(\ref{1.3})
is analyzed. The algorithm used here is a suitable modification of the algorithms in \cite{DC, LT} and it is 
composed of the following three steps:
\begin{itemize}
 \item [$\bullet$] {\bf Step 1}: solve a nonlinear problem over a coarse mesh with mesh size $H$ which provides an approximate solution, say  $\bu_H.$ 
 \item[$\bullet$] {\bf Step 2}: linearize the nonlinear system around the coarse grid solution $\bu_H$ and solve the resulting linearized problem over a fine mesh with mesh size $h$ and denote its solution as $\bu_h^{*}.$
 \item[$\bullet$]{\bf Step 3}: correct the  solution $\bu_h^{*}$ obtained in {\bf Step 2} over fine mesh which provides  an updated  final solution $\bu_h.$ 
\end{itemize}
As a result of the above mentioned three steps algorithm, the error $\|\bu-\bu_{h}\|$ is 
of the same order as $\|\bu-\tilde{\bu}_h\|$, where $\tilde{\bu}_h$ is the solution of 
the standard Galerkin system on a fine mesh $h$ with an appropriate scaling between $h$ and $H$. 

The two grid method has been extensively studied for Navier-Stokes equations by Layton \cite{LY}, Layton and Tobiska \cite{LT}, 
Layton and Lenferink \cite{LF1}-\cite{LF2}, Girault and Lions \cite{GL1,GL2}, 
Dai {\it et al.} \cite{DC}, Abboud {\it et al.} \cite{AGS1}-\cite{AGS}, Frutos {\it et al.} \cite{FGN}. 

In \cite{LT}, Layton {\it et al.} have examined a coarse mesh correction in the third  step for a steady state Navier-Stokes equations. 
But, this correction fails to improve the results obtained in {\bf Step 2} and as a result, 
optimal error estimate in $L^2$-norm for the velocity is obtained when $h=\mathcal{O}(H^{3/2}).$  
Based on stream function formulation, a two-grid finite element method has been studied by Fairag \cite{F}. All the above results 
have been discussed for the steady state Navier-Stokes equations on a convex polyhedra or on a convex polygon. 
Subsequently, Girault {\it et al} \cite{GL1} in their work on steady state 
Navier-Stokes equations have analyzed a two level two-grid algorithm and 
have obtained optimal $\bH^1$-norm error estimate for the velocity vector with a choice $h=\mathcal{O}(H^{2}),$ 
when the problem is defined on a Lipschitz polyhedron or on a convex polyhedron. The analysis is further extended to the 
 transient Navier-Stokes equations in \cite{GL2}, and optimal error estimate in $L^{\infty}(\bH^1)$-norm is 
established with a choice $h=\mathcal{O}(H^{2}),$  when $\Omega$ is a Lipschitz  polyhedron or a convex polyhedron. 
In both of these articles, the key approach is to 
exploit the contribution of the coarse grid solution in $L^3 (\Omega)$-norm.

In the context of nonlinear Galerkin method, two grid method is applied to the 2D-transient Navier-Stokes equations 
by Ait Ou Amni and Marion in \cite{AOM}. They have shown that the nonlinear Galerkin solution has the same accuracy as that of 
the standard Galerkin solution, both for velocity in $H^1$-norm and for pressure in $L^2$-norm with a choice $h=\mathcal{O}(H^{2})$. 
Further, they have 
penalized their two-grid algorithm to get rid of the coupling between velocity and pressure with penalization parameter $\epsilon$ and 
have recovered the same accuracy for the penalized two-grid Galerkin solution as that of the standard Galerkin solution with 
$h=\mathcal{O}(H^{2})$ and $h=\mathcal{O}(\epsilon^{1/2}).$

Garc\'{i}a-Archilla and Titi in \cite{AT} have applied Post-Processed method to the semilinear scaler elliptic equations in any 
dimensions and have derived optimal error bounds in $\bH^1$-norm for the post-processed solution with a choice 
$h=\mathcal{O}(H^{r+1}|log(H)|)^{1/r}))$, where the post-processed solution is approximated by the polynomials of degree $r$ with $r\geq 2$.       

Recently, Frutos {\it et al} \cite{FGN} have applied the two-grid scheme to the
incompressible Navier-Stokes equations using mixed-finite elements, the mini-element, 
the quadratic and the cubic Hood-Taylor elements for spatial discretization and a backward Euler method and 
a two step backward difference scheme for time discretization and have derived the rate of convergence 
of the fine mesh in the $\bH^1$-norm by taking $h=\mathcal{O}(H^2),$ which is an improvement over $ h=\mathcal{O}(H^{3/2})$ obtained in \cite{AGS}. 

In \cite{YH}, a fully discrete two-level method consisting of Crank-Nicolson extrapolation method with solution 
$(\bu_{H,\tau_0},p_{H,\tau_0})$ on a space-time coarse grid $J_{H,\tau_0}$ and a backward Euler method with solution 
$(\bu_{h,\tau},p_{h,\tau})$ on a space-time fine grid $J_{h,\tau}$ is discussed. They have obtained  convergence rate for the 
two level solution $(\bu_{h,\tau},p_{h,\tau}),$ which is of same order as that of the one level standard Crank-Nicolson extrapolation solution 
if $\tau_0^{3/2}+H^{3/2}=\mathcal{O}(\tau)$ for $t\in[0,1]$ and $\tau_0^{2}+H^{2}=\mathcal{O}(\tau)$ for $t\in[1,T]$.  

An attempt has been made in this article to discuss optimal error estimates in $L^{\infty}({\bL^2})$ and $L^{\infty}(\bH^1)$-norms 
for the velocity and $L^{\infty}(L^2)$-norm for the pressure using a three level two-grid finite element 
method for the 2D-transient Navier-Stokes equations. The major contributions are given in terms of the following two tables. 
In Table 1, we present the order of convergences for the two-grid algorithm (\ref{1a1})-(\ref{3a1}) stated in Section {\bf 3} 
for the pair of finite element spaces $(\bH_\mu , L_\mu ), \mu = H, h$ satisfying the approximation properties mentioned 
in {\bf (B1)}-{\bf (B2)}. Table 2 provides the largest scaling between coarse and fine meshes 
for which the desired fine mesh accuracy is obtained for both velocity and pressure. 
 {{
\begin{table}[ht!]
 \centering
\begin{tabular}{|c|c|c|c|c|c|c|c|c|}
 \hline
  Solution   &{Velocity } & {Velocity } & {Pressure}  \\
             &    {in $\bL^2$-norm}                            &    { in $\bH^1$-norm}      &   {in $L^2$-norm} \\
 \hline
 \hline  $(\bu-\bu_H,p-p_H)$ &  $H^2$  & $H$          & $H$  \\
 \hline $(\bu-\bu_h^{*},p-p_h^{*})$ &\qquad $h^2+H^{3-\delta}$ &\qquad $h+H^{3-\delta}$&\qquad$h+H^{3-\delta}$\\
 \hline  $(\bu-\bu_h,p-p_h)$ & \qquad$h^2+H^{4-2\delta}$&\qquad$h+H^{4-\delta}$&\qquad $h+H^{4-\delta}$\\
 \hline
 \end{tabular}
 \caption{ Error estimates obtained from the two-grid algorithm for arbitrarily small $\delta>0$.}
 \end{table}
 }
 {{
\begin{table}[ht!]
 \centering
 \begin{tabular}{|c|c|c|c|c|c|c|c|c|}
 \hline
  Solution   &{Velocity } &{Velocity } & {Pressure}  \\
             &    {in $\bL^2$-norm}                            &  { in $\bH^1$-norm}      &  {in $L^2$-norm} \\
 \hline
 \hline  $(\bu-\bu_H,p-p_H)$ &$H^2$  &$H$          &$H$  \\
 \hline $(\bu-\bu_h^{*},p-p_h^{*})$ & \qquad$h\sim H^{(3-\delta)/2}$ & \qquad$h\sim H^{3-\delta}$&\qquad$h\sim H^{3-\delta}$\\
 \hline  $(\bu-\bu_h,p-p_h)$ & $h\sim H^{2-\delta}$&\qquad$h\sim H^{4-\delta}$& \qquad$h\sim H^{4-\delta}$\\
 \hline
 \end{tabular}
 \vspace{.1cm}
 \caption{The largest scaling for optimal error estimates.}
 \end{table}
 }}
 \newline
 It is observed  from Tables 1 and 2 that the introduction of Step {\bf 3} leads to a good improvement in scaling 
between $H$ and $h$ for both $\bH^1$-norm for the velocity and $L^2$-norm for the pressure, that is, the scaling improves from 
$h\sim H^{3-\delta}$ to $h\sim H^{4-\delta}$. It also improves the scaling for 
$\bL^2$-norm of the velocity from Step {\bf 2} to Step {\bf 3} from $h\sim H^{(3-\delta)/2}$ to $h\sim H^{(2-\delta)}$ for arbitrarily small $\delta>0.$ 

The main contributions of this paper can be summarized as follows:
{\begin{itemize}
\item [(i)] Based on the steady state Oseen projection and Sobolev estimates in Lemma \ref{nlt} involving $\delta$, 
optimal error estimates for the two-grid Galerkin approximations to the velocity 
in $L^{\infty}({\bf H}^1)$-norm and to the pressure in $L^{\infty}( L^2)$-norm 
with the largest scaling between $H$ and $h$, $h\sim H^{3-\delta}$ and $h\sim H^{4-\delta}$ 
for {\bf Step 2} and {\bf Step 3}, respectively are derived. The result obtained in {\bf Step 2} is an improvement over the result obtained by Frutos {\it et al} \cite{FGN}. They have obtained using 
first order mini-elements the largest scaling between $H$ and $h$,  as $h\sim H^2$ for both 
$L^{\infty}({\bf H}^1)$-norm for the velocity and $L^{\infty}( L^2)$-norm for the pressure. 
\item[(ii)] A use of linearized backward Oseen problem with related estimates yields optimal $L^{\infty}({\bf L}^2)$-norm estimates for the velocity in {\bf Step 2} with a choice 
$h=\mathcal{O}(H^{(3-\delta)/2})$ and {\bf Step 3}  with a choice $h=\mathcal{O}(H^{2-\delta})$ for $\delta>0$ arbitrarily small.
\item[(iii)] Under the assumption of uniqueness condition,  {\it a priori} error estimates are obtained which hold uniformly in time.
\end{itemize}}

The remaining part of the paper consists of the following sections. In Section {\bf 2}, some preliminaries to be 
used in the subsequent sections are presented. In Section {\bf 3}, semidiscrete two-grid 
finite element approximations are introduced. Optimal error estimates for velocity and pressure are established 
in Section {\bf 4}. 
Section {\bf 5} deals with the backward Euler method applied 
to the semidiscrete two grid
system. Finally, in Section {\bf 6}, the results of some numerical examples which confirm our theoretical results are presented. 
\section{Preliminaries}
We denote $\mathbb {R}^2$-valued function spaces using bold 
face letters, that is, ${\bf H}_0^1 = (H_0^1(\Omega))^2$, ${\bf L}^2 = (L^2(\Omega))^2$ and ${\bf H}^m=(H^m(\Omega))^2.$ 
The standard notations for Lebesgue and Sobolev spaces with their norms are employed in the paper. 
The space $\bH_0^1$ is equipped with a norm $\|\nabla\bv\|= \left({\sum_{i,j=1}^{2}}
(\partial_j v_i, \partial_j v_i)\right)^{1/2}$.
Given a Banach space $X$ endowed with norm ${\parallel \cdot \parallel}_X$, 
let $L^p(0,T;X)$ be the space of all strongly measurable functions $\phi:[0,T] \rightarrow X$ 
satisfying $\displaystyle{\int_0^{T}} {\parallel \phi(s) \parallel}^p_X\,ds< \infty$ and 
for $p=\infty$, $\displaystyle{\mathop{ess \sup}_{t \in[0,T]} {\parallel \phi(t) \parallel}_X< \infty}$. 
Also, define
 \begin{eqnarray}
 {\bf J}&=& \{\bphi \in {\bf {L}}^2 :\nabla \cdot \bphi = 0\;\;
{\mbox {\rm in}}
 \;\; \Omega, \;\;\bphi \cdot {\bf {n}} |_{\partial \Omega} = 0\;\;
 {\mbox {\rm holds} \;  {\rm weakly}} \}\nonumber\\
{\bf J}_1 &=& \{{\bphi} \in {\bf {H}}_0^1 : \nabla \cdot \bphi = 0\},\nonumber
 \end{eqnarray}
where ${\bf {n}}$ is the unit outward normal to the boundary
$\pr \Omega$ and $\bphi \cdot {\bf {n}} |_{\pr \Omega} = 0$ should
be understood in the sense of trace in $\bH^{-1/2}(\partial \Omega)$,
see \cite{temam}. Let $H^m/{\rm I\!R} $ be the quotient space with norm $\| \phi\|_{H^m /{\rm I\!R}}
 = \inf_{c\in{\rm I\!R} }\| \phi+c\|_m$. For $m=0$, it is denoted by $L^2 /{\rm I\!R}$.\\
\noindent
Throughout this paper, we make the following assumptions: \\
\noindent 
({\bf A1}). For ${\bf {g}} \in \bL^2$, let $(\bv, q) \in \bJ_1\times L^2 /{\rm I\!R}$ be the unique solution 
to the steady state Stokes problem
\be
&&-\Delta {\bv} + \nabla q = {\bf {g}},\\
&&\nabla \cdot\bv = 0\;\;\; {\mbox {\rm in} }\;\;\; \Omega,
\;\;\;\; \bv|_{\pr \Omega} = 0
\ee
 satisfying the regularity result \cite{temam}:
\begin{align}
\| \bv \|_2 + \|q\|_{H^1 /{\rm I\!R}} \le C\|{\bf {g}}\|.\nonumber
\end{align}
It is easy to show that 
\ben \label{2.1*} 
\|\bv\|^2 \le \lambda_1^{-1} \| \nabla \bv \|^2
\;\; \forall \bv \in {\bf H}_0^1(\Omega),
\een
where $ \lambda_1 $ is the minimum eigenvalue of the Laplacian with zero Dirichlet boundary  condition. 

\noindent
({\bf A2}). There exists a positive constant $M_0$ such that the initial velocity $\bu_0$ and external force $\f$ satisfy for 
$t\in (0,T]$ with $0<T< \infty$ 
\begin{eqnarray}
\bu_0 \in \bJ_1\cap\bH^2,\,\, \f,\,\f_t \in L^{\infty}(0,\,T;\,\bL^2)\,\, \rm{with}\,\, \|\bu_0\|_2 \le M_0,\,\, 
\displaystyle{\mathop{ess\,\sup}_{0<t\leq T}}\,\{\|\f(\cdot,t)\|,~\|\f_t(\cdot,t)\|\}\leq M_0.\nonumber
\end{eqnarray}
\noindent
Now, for $\bv, \bw, \bphi \in \bH_0^1$, define
$
a(\bv, \bphi) := (\nabla \bv, \nabla \bphi)
$
and
$
b(\bv, \bw,\bphi):= \frac{1}{2} (\bv \cdot \nabla \bw , \bphi)
- \frac{1}{2} (\bv \cdot \nabla \bphi, \bw).
$\\
\noindent
The weak formulation of (\ref{1.1})-(\ref{1.3}) is to find $(
\bu(t), p(t))\in \bH_0^1\times L^2 /{\rm I\!R}$, such that $\bu(0)= \bu_0$ and for $t>0$ 
\begin{eqnarray} \label{2.2atw}
 \left.
 \begin{array}{rcl}
&&(\bu_t, \bphi) +\nu\, a( \bu, \bphi)+b( \bu,\bu, \bphi)-( p,\nabla \cdot \bphi) 
=({\bf f},\bphi)
\;\;\;\;\; \forall \bphi \in {\bf H}_0^1,\\
&&(\nabla \cdot \bu,\chi)=0\;\;\;\forall \chi\in L^2.
 \end{array}
 \right\}
 \end{eqnarray}
\noindent
Equivalently, find  $\bu(t) \in {\bf J}_1 $ such that for $\bu(0)= \bu_0$, $t>0$,
\begin{align}
\label{2.3atw}
(\bu_t, \bphi) +\nu\,a( \bu, \bphi )+b( \bu,\bu, \bphi)= ({\bf f},\bphi)\;\;\;\forall \bphi \in {\bf J}_1.
 \end{align}
We recall below, the following regularity results.
\begin{ldf}\label{L45}\cite[pp. 285, 302]{HR82}
Let the assumptions {\rm{({\bf A1})-({\bf A2})}} hold true. Then, for any $T$ with $0<T< \infty$, for any 
fixed $\alpha>0$ and for some constant $C=C(M_0)$, the solution of 
(\ref{2.2atw}) satisfies
\begin{align}
& \displaystyle{\sup_{0<t\leq T}} \{\|\bu(t)\|_2+\|\bu_t(t)\|+\|p(t)\|_{\bH^1/{\mathbb R}}\}\leq C,\nonumber\\
&\bs^{-1}(t)\displaystyle{\int_0^t}e^{2\alpha \tau}(\| \bu(\tau)\|_2^2+\|p(\tau)\|_1^2)d\tau\leq C,\nonumber\\
 & \displaystyle{\sup_{0<t\leq T}} e^{-2\alpha t}\int_0^t e^{2\alpha\tau}\|\bu_t(\tau)\|_1^2d\tau\leq C, \,\,\, 
\displaystyle{\sup_{0<t\leq T}} \tau(t)\|\bu_t(t)\|_1^2\leq C,\nonumber\\
&\displaystyle{\sup_{0<t\leq T}}e^{-2\alpha t}\displaystyle{\int_0^t}\bs(\tau)(\|\bu_t(\tau)\|_2^2+\|\bu_{\tau\tau}(\tau)\|^2
+\|p_\tau(\tau)\|^2_{\bH^1/{\mathbb R}})d\tau\leq C,\nonumber
\end{align}
where $\tau(t):=\min\{t,\,1\}$ and $\bs(t):=\tau(t)e^{2\alpha t}$.
\end{ldf}
\section{Two-Grid Formulation}
\setcounter{equation}{0}

Consider two admissible shape regular finite triangulations of $\bar\Omega\,$: a coarse mesh $\mathcal{T}_H$ with mesh 
size $H$ and a fine mesh $\mathcal{T}_h$ with mesh size $h$, where $h\ll H$. Let ${\bf H}_\mu$ and $L_\mu$ be the 
finite dimensional subspaces of ${\bf H}_0^1 $ and $L^2$, respectively, where $\mu=H, h$. 
Let us also consider the associated divergence free subspaces 
$\bJ_\mu$ of $\bH_\mu$,  where
${\bf J}_\mu = \{ \bphi_\mu \in {\bf H}_\mu: (\nabla \cdot \bphi_\mu,\chi_\mu)
=0 \;\;\; \forall \chi_\mu\in L_\mu \}.$
Note that $\bJ_\mu$ is not a  subspace of $\bJ_1$.

Let the spaces ${\bf H}_\mu$ and $L_\mu$ satisfy the following  properties:\\
\noindent
({\bf B1}). ({\it Approximation property}) For $\bw \in {\bf H}_0^1 \cap {\bf {H}}^2 $
and $ q \in H^1/ {\rm I\!R}$, there exist approximations $i_\mu \bw \in {\bH}_\mu$ and $ j_\mu q \in L_\mu$, such that
\be
\|\bw-i_\mu\bw\|+ \mu \| \nabla (\bw-i_\mu \bw)\| \le C \mu^2
\| \bw\|_2, \;\;\;\; \| q - j_\mu q
 \|_{L^2 /{\rm I\!R}} \le C \mu \| q\|_{H^1 / {\rm I\!R}}.
\ee
\noindent
({\bf B2}). ({\it Uniform inf-sup condition}) There exists a positive constant $C$, independent of $\mu$, such that 
\be
 \displaystyle{\sup_{{\bphi}_\mu\in\bH_\mu\diagdown\{0\}}}\frac{|(q_\mu, \nab \cdot {\bphi}_\mu)|}{\|\nab {\bphi}_\mu \|} \ge C \| q_\mu\|_{
L_\mu/N_\mu}\,\,\forall q_\mu \in L_\mu,
\ee
where
$N_\mu= \{ q_\mu\in L_\mu:\forall \bphi_\mu\in {\bf H}_\mu, (q_\mu, \nab \cdot {\bphi}_\mu ) = 0\}$.

Note that $\bJ_\mu$ is not a  subspace of $\bJ_1$. With $P:\bL^2\rightarrow \bJ$ an orthogonal projection, 
set the Stokes operator 
$\tilde\Delta=P\Delta$. The $L^2$ projection $P_\mu:\bL^2\rightarrow \bJ_\mu$ 
satisfies the following properties \cite{HR82}: 
\begin{eqnarray} \label{3.4}
 \left.
 \begin{array}{rcl}
&&\|\bphi- P_\mu \bphi\|+ \mu \|\nabla  P_\mu \bphi\| \leq C \mu
 \|\nabla \bphi\|\qquad\qquad \bphi \in \bJ_1,\\
&&\|\bphi- P_\mu \bphi\|+ \mu\|\nabla  (\bphi-P_\mu \bphi)\|\leq C\mu^2\|\tilde\Delta\bphi\|\,\,\,\, \bphi\in\bJ_1\cap\bH^2.
 \end{array}
 \right\}
 \end{eqnarray}
Define the discrete analogue of the Stokes operator as $\tilde \Delta_\mu=P_\mu\Delta_\mu$, 
where $\Delta_\mu$ is defined by $(\Delta_\mu \bv_\mu,\bphi_\mu)=-(\nabla\bv_\mu,\nabla\bphi_\mu)$, for all 
$\bv_\mu$, $\bphi_\mu\in \bH_\mu$. Define the 'discrete' 
Sobolev norms on $\bJ_\mu$ (see \cite{HR82}) as for $r\in {\mathbb R}$ and 
for $\bv_\mu\in\bJ_\mu$, $\|\bv_\mu\|_r:=\|(-\tilde\Delta_\mu)^{r/2}\bv_\mu\|$.

\noindent 
The operator $b(\cdot, \cdot, \cdot)$ satisfies the
antisymmetric property; that is,
\begin{eqnarray}\label{3.1*}
b(\bv_\mu, \bw_\mu, \bw_\mu) = 0 \;\;\; \forall \bv_\mu, \bw_\mu \in {\bH}_\mu.
\end{eqnarray}
In the following lemma, we state without proof some estimates of the trilinear term $b(:,.,.)$. For a proof, 
see \cite[pp 360]{HR90} and \cite[pp. 2044]{LT}.
\begin{ldf}\label{nlt}
The trilinear form $b(\cdot,\cdot,\cdot)$ satisfies the following estimates:
\begin{eqnarray*} 
|b(\bphi,\bxi,\chi)|\leq C
 \left\{
 \begin{array}{rcl}
&&\hspace{-.8cm} \| \nabla {\bphi}\|^{1/2}\|\tl{\Delta}_\mu {\bphi}\|^{1/2} \| \nabla  \bxi\|\;\| \chi\|,~\rm{for~all}\,{\bphi},
~ \bxi,~ \chi \in  {\bH}_\mu,\\
&&\hspace{-.8cm} \| \nabla \bphi\|\|\nabla \bxi\|^{1/2}
\|\tl{\Delta}_\mu {\bxi}\|^{1/2} \| \chi\|,~\rm{for~all}\,{\bphi},~ \bxi,~ \chi \in  {\bH}_\mu,\\
 &&\hspace{-.8cm}\|\bphi\|\|\nabla \bxi\|
\|\tl\Delta \chi\|,~\rm{for\,all}~  {\bphi},~ \bxi\in {\bH}_0^1,~ \chi \in {\bH}_0^1\cap{\bH}^2,\\
&&\hspace{-.8cm}\|\nabla\bphi\|\|\bxi\|
\|\tl\Delta \chi\|,~\rm{for~all}~  {\bphi},~ \bxi\in {\bH}_0^1,~ \chi \in {\bH}_0^1\cap{\bH}^2,\\
 && \hspace{-.8cm}\|\nabla\bphi\|\|\nabla\bxi\|\|\nabla \chi\|,~\rm{for~all}~ {\bphi},\,\bxi,\,\chi\in {\bH}_0^1,\\
 &&\hspace{-.8cm}\|\bphi\|^{1-\delta}\|\nabla \bphi\|^{\delta}\|\nabla\bxi\|\|\nabla \chi\|,~ 
\rm{for~all}~ {\bphi},\,\bxi,\,\chi\in {\bH}_0^1,
 \end{array}
 \right.
 \end{eqnarray*}
where $\delta>0$ is arbitrarily small.
\end{ldf}
The three level two-grid  semidiscrete algorithm applied to (\ref{1.1})-(\ref{1.3}) is described as follows: \\
\noindent
{\it {\bf Step 1} ( nonlinear system (\ref{1.1}) on a coarse grid): Find $\bu_H\in \bJ_H$ such 
that for all $\bphi_H\in \bJ_H$ for $\bu_H(0)=P_H\bu_{0}$ and $t>0$
\begin{eqnarray}\label{1a1}
(\bu_{Ht},\bphi_H)+\nu a(\bu_H, \bphi_H) +b(\bu_H,\bu_H,\bphi_H)=({\bf f},\bphi_H).
 \end{eqnarray}
{\bf Step 2} ( Update on a finer mesh with one Newton iteration ) : Seek $\bu_{h}^{*}\in \bJ_h$ 
such that for all $\bphi_h\in \bJ_h$ for $\bu_h^{*}(0)=P_h\bu_{0}$ and $t>0$
 \begin{align}\label{2a1}
(\bu_{ht}^{*},\bphi_h)+\nu a(\bu_h^{*}, \bphi_h) 
+b(\bu_h^{*},\bu_H,\bphi_h)+b(\bu_H,\bu_h^{*},\bphi_h)
=({\bf f},\bphi_h)+b(\bu_H,\bu_H,\bphi_h).
\end{align}
{\bf Step 3} ( Correction on a fine mesh) : Find $\bu_{h}\in \bJ_h$ such that for all $\bphi_h\in \bJ_h$ 
for $\bu_h(0)=P_h\bu_{0}$ and $t>0$	
\begin{align}\label{3a1}
(\bu_{ht},\bphi_h)&+\nu a(\bu_h, \bphi_h) +b(\bu_h,\bu_H,\bphi_h)
+b(\bu_H,\bu_h,\bphi_h)\nonumber\\
&=({\bf f},\bphi_h)+b(\bu_H,\bu_h^{*},\bphi_h)+b(\bu_h^{*},\bu_H-\bu_h^{*},\bphi_h).
\end{align}
}
The following inequality  will be used frequently in our error analysis:
\begin{align}\label{tw11}
 \inf_{\bphi_h\in\bJ_h}\sup_{\bv_h\in \bJ_h}
\frac{\nu\,a(\bphi_h,\bv_h)+b(\bu_H,\bphi_h,\bv_h)+b(\bphi_h,\bu_H,\bv_h)}{\|\nabla\bphi_h\|\|\nabla\bv_h\|}\geq\gamma>0.
\end{align}
For a proof, see \cite{LT}.

For uniform estimates in time, we shall further assume  the following 
uniqueness condition: 
\begin{align}\label{tn}
 \frac{N}{\nu^2}\|\f\|_{L^{\infty}(0,\infty;\bL^2)}<1\,\,\, \text{and}\,\,N=
 \sup_{\bu,\bv,\bw\in \bH_0^1(\Omega)}\frac{b(\bu,\bv,\bw)}{\|\nabla\bu\|\|\nabla\bv\|\|\nabla\bw\|},
\end{align}

The main results of this section are stated in the following theorems.
\begin{tdf}\label{T5}
 Let $\Omega$ be a convex polygon and let assumptions {\rm (\bf {A1})-(\bf{A2})} and {\rm (\bf{B1})}-{\rm(\bf{B2})} hold true. 
Further, let the discrete initial velocity $\bu_{0h}\in \bJ_h$ with $\bu_{0h}=P_h\bu_0$. 
Then, there exists a positive constant $C$, independent of $h$, such that for $t\in (0,T]$ with $0<T<\infty$, 
the following estimates hold true:
\begin{align} 
 \label{error-velocity}
\| (\bu-\bu_h)(t)\|\leq K(t)(h^2+ H^{4-2\delta}\big),\,\,
\|\nabla(\bu-\bu_h)(t)\|\leq K(t) (h+H^{4-\delta}),
\end{align}
and
\begin{align} 
\label{error-pressure}
\| (p-p_h)(t)\|\leq K(t)(h + H^{4-\delta}\big),
\end{align}
where $\delta>0$ is arbitrarily small and $K(t)=C e^{Ct}$. Under uniqueness condition (\ref{tn}), $K(t)=C$ and 
the estimates in Theorem \ref{T5} are valid uniformly in time.
\end{tdf}
The remaining part of this paper is devoted to the derivation of results, which will lead to the proof of Theorem 
\ref{T5}.
\section{Error Estimates}
This section deals with  optimal error estimates of the semidiscrete two-grid algorithm.
Since $\bJ_h$ is not a subspace of $\bJ_1$, the weak solution $\bu$ satisfies
\begin{eqnarray}\label{4a1}
(\bu_t,{\bphi}_h)+\nu\, a(\bu,{\bphi}_h)+b(\bu,\bu,{\bphi}_h)=(\f,\bphi_h)+(p,\nabla \cdot {\bphi}_h) \,\,\,\, \forall {\bphi}_h \in \bJ_h.  
\end{eqnarray}
Define $\e_H :=\bu-\bu_H$, $\e^{*}:=\bu-\bu_h^{*}$ and $\e_h:=\bu-\bu_h$. 
Then, a use of  (\ref{2a1}) and (\ref{4a1}) yields
\begin{align}\label{5}
(\e_t^{*},{\bphi}_h)&+\nu\, a(\e^{*},{\bphi}_h) +b(\e^{*},\bu_H,\bphi_h)+b(\bu_H,\e^{*},\bphi_h)\nonumber\\
&=-b(\e_H,\e_H,\bphi_h)+(p,\nabla \cdot {\bphi}_h)\,\,\,\,\,\,\,\,\qquad \forall {\bphi}_h \in \bJ_h.
\end{align}
Subtract (\ref{3a1}) from (\ref{2a1}) and then add the resulting equation to (\ref{5}) to arrive at
 \begin{align}\label{7}
(\e_{ht},{\bphi}_h)&+\nu\, a(\e_h,{\bphi}_h)+b(\e_h,\bu_H,\bphi_h)+b(\bu_H,\e_h,\bphi_h)=-b(\e_H,\e^{*}, \bphi_h)\nonumber\\
&-b(\e^{*},\e_H, \bphi_h)
+b(\e^{*},\e^{*},\bphi_h)+(p,\nabla\cdot \bphi_h) \,\,\,\,\,\,\,\, \forall {\bphi}_h \in \bJ_h.
\end{align}
For analyzing optimal error estimates of $\e_h$ in $L^{\infty}(\bL^2)$ and $L^{\infty}(\bH^1)$-norms, 
define  an auxiliary projection $\tilde\bu_h(t)\in \bJ_h$, $0<t\leq T,$ for a given $\bu,$ as a solution of the following modified steady state Oseen  problem: 
\begin{eqnarray}\label{E616*}
\nu\, a(\bu-\tilde\bu_h,{\bphi}_h)+b(\bu_H,\bu-\tilde\bu_h,\bphi_h)+b(\bu-\tilde\bu_h,\bu_H,\bphi_h)
 =(p,\nabla\cdot \bphi_h)
\,\,\,\, \rm{for~ all} ~ {\bphi}_h\in \bJ_h.
\end{eqnarray}
Now split $\e_h$ as
\begin{eqnarray}\label{17}
\e_h=:(\bu-\tilde\bu_h)+(\tilde\bu_h-\bu_h):=\bzeta+\Th,
\end{eqnarray}
where $\bzeta:=\bu-\tilde\bu_h$ and $\Th:=\tilde\bu_h-\bu_h$.\\
\noindent
A use of (\ref{7})-(\ref{17}) leads to  
\begin{align}\label{B18*}
 (\Th_t,\bphi_h)&+\nu\, a(\Th,\bphi_h)+b(\Th,\bu_H,\bphi_h)+b(\bu_H,\Th, \bphi_h)=-(\bzeta_{t},\bphi_h)\nonumber\\
 &-b(\e_H,\e^{*},\bphi_h)-b(\e^{*},\e_H,\bphi_h)+b(\e^{*},\e^{*},\bphi_h).
\end{align}
To seek estimates for $\Th$, we need  estimates for $\bzeta$, $\e_H$ and $\e^{*}$, which appear on the right hand side of (\ref{B18*}).
\begin{ldf}\label{1l1*}\cite[see page 362, proposition 3.2]{HR90}
Let $\bu_H(t)$ be the solution of (\ref{1a1}) in $[0,T),\,\,0<T<\infty$ satisfying 
$\bu_{0H}=P_H\bu_{0}$, and let the assumptions {\rm{\bf (A1)}}-{\rm{\bf (A2)}} hold true. 
Then, there exists a positive constant $C=C(\gamma,\nu,\alpha,\lambda_1, M_0)$ such that  
the following hold true for all $t>0$
\begin{align}
\|\bu_H(t)\|_2+\|\bu_{Ht}(t)\|+e^{-2\alpha t}\displaystyle{\int_0^t} e^{2\alpha s}(\|\tilde\Delta \bu_H(s)\|^2+
\|\nabla\bu_{Ht}(s)\|^2) ds
\leq C,\nonumber 
\displaystyle{\sup_{0<t<T}} \tau(t)\|\bu_{Ht}(t)\|_1^2\leq C.
\end{align}
 \end{ldf}
\begin{ldf}\label{T2}\cite[\it estimates for $\e_H$]{HR82}
Let the assumptions {\rm (\bf {A1})-(\bf{A2})} and {\rm (\bf{B1})}-{\rm(\bf{B2})} hold true. With initial velocity 
$\bu_{0H}=P_H\bu_0$, let the discrete solution pair $(\bu_H(t),p_H(t))$ satisfies (\ref{1a1}). Then, there exists a 
positive constant $C$, independent of $H$, such that for $0<t\leq T$
 \begin{align}
 \|\e_H(t)\|\leq K(t)\,H^2 \,\,\,\,\qquad\|(p-p_H)(t)\|\leq K(t)\,H\nonumber
\end{align}
and 
\begin{align}
\bs^{-1}(t)\int_0^t e^{2\alpha \tau}\|\e_H(\tau)\|^2 d\tau\leq K(t)\,H^4,\,\,
e^{-2\alpha t}\int_0^t \bs(\tau)\|\nabla\e_{H\tau}(\tau)\|^2 d\tau\leq 
K(t)\,H^2,\nonumber
\end{align}
where $K(t)=C e^{Ct}$. If in addition, uniqueness condition (\ref{tn}) holds true, then $K(t)=C$ and the results are 
valid uniformly in time.
\end{ldf}
In order to derive estimates for $\e^{*}$, split it as
\begin{align}\label{tw9*}
\e^{*}=(\bu-\tilde\bu_h)+(\tilde\bu_h-\bu_h^{*}):=\bzeta+\brho.
\end{align}
The following lemma provides estimates for $\tilde\bu_h.$  Since a suitable modification of the proofs in \cite{LT} will provide 
a proof, we state below the results without the proofs.
\begin{ldf}\label{nl2}
Let the approximation property {\rm{\bf (B1)}} be satisfied. Further, let $\bu$ and $\tilde\bu_h$ be the solution of (\ref{2.2atw}) and 
(\ref{E616*}), respectively. Then, the following estimates hold true:
\begin{align}
 &\|\nabla\bzeta(t)\| \leq C h\,{\mathcal K}(t),\,\,\|\nabla \bzeta_t(t)\|\leq C h \,
 \left(\mathcal{K}_t(t)+{\mathcal K}(t)\|\nabla\bu_{Ht}\|\right),\nonumber\\
 &\|\bzeta(t)\|\leq C h{\mathcal K}(t)(h+\|\e_H\|),\,\,\, 
 \|\bzeta_t(t)\|\leq Ch(h+\|\e_H\|)({\mathcal K}(t)\|\nabla\bu_{Ht}\|+{\mathcal K}_t(t)),\nonumber
\end{align}
where ${\mathcal K}(t):=\|\tilde\Delta\bu(t)\|+\|\nabla p(t)\|$ and 
${\mathcal K}_t(t):=\|\tilde\Delta\bu_t(t)\|+\|\nabla p_t(t)\|$. 
\end{ldf}

We note that ${\mathcal K}(t)$ and ${\mathcal K}_t(t)$ satisfy  for some positive constant $K$ 
\begin{eqnarray} \label {bound-K}
\displaystyle{\sup_{0 < t\leq T}}\;\; e^{-2\alpha t} \;\int_{0}^t e^{2\alpha s} \Big({\mathcal K}^2(s) + \tau (s){\mathcal K}_t^2(s)\Big)\; ds \;\leq  K.
\end{eqnarray}

Below in Subsections {\bf 4.1} and {\bf 4.2}, we focus on  the semidiscrete error estimates related to {\bf Step 2} and {\bf Step 3}.
\subsection{Error estimates for Step 2}
In this subsection, the semidiscrete error estimates corresponding to {\bf Step 2} are derived.
\begin{ldf}\label{nl5}
With $\displaystyle{0\leq\alpha \leq \gamma\lambda_1}$ and $\beta= \gamma-\alpha\lambda_1^{-1}>0,$ let the hypothesis of Lemma \ref{T2} hold true. Then, 
the following estimate holds:
\begin{align}
\beta \;\bs^{-1}(t)\int_0^t e^{2\alpha \tau} \|\nabla \e^{*}(\tau)\|^2 d\tau \leq\;  K(t) (h^2+H^{6-2\delta}).\nonumber
\end{align} 
Here and elsewhere in this paper, $K(t)$ denotes $Ce^{Ct}$ and under uniqueness assumption (\ref{tn}), $K(t)$, reduces to a positive constant $K.$
\end{ldf}
 \noindent
{\it Proof.}  Choose $\bphi_h=P_h \e^{*}=(P_h \bu-\bu)+\e^{*}$ in (\ref{5}). Then, a use of (\ref{tw11}) 
yields
\begin{align}\label{ne21}
\frac{1}{2}\frac{d}{dt}\|\e^{*}&\|^2+\gamma\,\|\nabla\e^{*}\|^2\leq (\e^{*}_t,\bu-P_h \bu)
+\nu\,a(\e^{*},\bu-P_h\bu)+b(\bu_H,\e^{*},\bu-P_h\bu)\nonumber\\
&+b(\e^{*},\bu_H,\bu-P_h\bu)-b(\e_H,\e_H,P_h \e^{*})+(p,\nabla\cdot P_h\e^{*}).
\end{align}  
Using the definition of $P_h$, the first term on the right hand side of (\ref{ne21}) can be treated as
\begin{align}\label{ne22}
(\e^{*}_t,\bu-P_h \bu) &=(\bu_t-P_h \bu_t+P_h \bu_t-\bu_h^{*},\bu-P_h\bu)\nonumber\\
&=(\bu_t-P_h\bu_t,\bu-P_h \bu)=\frac{1}{2}\frac{d}{dt}\|\bu-P_h\bu\|^2.
\end{align}
An application of the Cauchy-Schwarz's inequality with (\ref{3.4}) leads to 
\begin{align}\label{ne24}
\nu|a(\e^{*},\bu-P_h\bu)|&\leq C(\nu)\|\nabla \e^{*}\|\|\nabla(\bu-P_h\bu)\|\leq C(\nu) h\|\tilde \Delta \bu\|\|\nabla\e^{*}\|.
\end{align}
Apply (\ref{3.4}), Lemma \ref{nlt} and the boundedness of $\|\nabla\bu_H\|$ to arrive at
\begin{align}\label{ne24*}
 |b(\bu_H,\e^{*},\bu-P_h\bu)+b(\e^{*},&\bu_H,\bu-P_h\bu)|\leq C\|\nabla \bu_H\|\|\nabla\e^{*}\|
\|\nabla(\bu-P_h\bu)\|\nonumber\\
&\leq Ch\|\tilde\Delta\bu\|\|\nabla \bu_H\|\|\nabla\e^{*}\|\leq Ch\|\tilde\Delta\bu\|\|\nabla\e^{*}\|.
\end{align}
The discrete incompressibility condition shows that
 \begin{align}\label{ne23*}
 |(p-j_h p,\nabla\cdot P_he^{*})|&\leq \|p-j_h p\|\|\nabla \e^{*}\|\leq C h\|\nabla p\|\|\nabla \e^{*}\|.
\end{align}
A use of (\ref{3.4}) with Lemma \ref{nlt} yields
\begin{align}\label{ne23}
 |b(\e_H,\e_H,P_h \e^{*})|&\leq C \|\e_H\|^{1-\delta}\|\nabla\e_H\|^{1+\delta}\|\nabla\e^{*}\|.
\end{align}
Substitute (\ref{ne22})-(\ref{ne23}) in (\ref{ne21}) along with the Young's inequality, (\ref{2.1*}) and 
multiply the resulting inequality by $e^{2\alpha t}$ to obtain
\begin{align}\label{ne25}
 \frac{1}{2}\frac{d}{dt}&e^{2\alpha t}\|\e^{*}\|^2+\frac{1}{2}\left(\gamma-\frac{\alpha}{\lambda_1} \right) e^{2\alpha t}\|\nabla \e^{*}\|^2
\leq \frac{1}{2} \frac{d}{dt}e^{2\alpha t}\|\bu-P_h\bu\|^2 \nonumber\\
&+C e^{2\alpha t}\left( h^2\mathcal{K}^2+\|\e_H\|^{2(1-\delta)}\|\nabla\e_H\|^{2(1+\delta)}\right)
-\alpha e^{2\alpha t}\|\bu-P_h\bu\|^2.
\end{align}
The last term in the right hand side of (\ref{ne25}) is negative. We drop this term. Integrate (\ref{ne25}) with respect to time, 
use $\|\e^{*}\|\leq \|\bu-P_h \bu\|$, $\|\e^{*}(0)\|=\|\bu_0-P_h \bu_0\|$ and Lemma \ref{L45} with $\beta=\gamma-\alpha\lambda_1^{-1}>0$ to arrive at  
  \begin{align}\label{ne25*}
 \beta  \displaystyle{\int_0^t}e^{2\alpha \tau}&\|\nabla \e^{*}(\tau)\|^2 d\tau\leq C\int_0^t e^{2\alpha \tau}\left(h^2\mathcal{K}^2(\tau)
+\|\e_H(\tau)\|^{2(1-\delta)}\|\nabla\e_H(\tau)\|^{2(1+\delta)}\right)d\tau\nonumber\\
&\leq C \left(h^2 \bs +\|\e_H\|_{L^{\infty}(\bL^2)}^{-2\delta}\|\nabla\e_H\|_{L^{\infty}(\bL^2)}^{2(1+\delta)}
\int_0^te^{2\alpha \tau}\|\e_H(\tau)\|^2d\tau\right).
 \end{align}
An application of Lemma \ref{T2} in (\ref{ne25*}) completes the proof.\hfill{$\Box$}

Next, we prove $L^{\infty}(\bH^1)$-norm estimate of $\e^{*}$.
\begin{ldf}\label{nl5*}
Under the hypotheses of Lemma \ref{nl5}, the following estimate holds true:
\begin{align}
 \bs^{-1}(t)\int_0^t\bs(\tau)\|\e^{*}_\tau(\tau)\|^2d\tau+\|\nabla \e^{*}\|^2\leq K(t)\,(h^2+H^{6-2\delta}).\nonumber
\end{align}
\end{ldf}
\noindent{\it Proof.} Substitute $\bphi_h=\bs P_h \e_t^{*}=\bs(P_h \bu_t-\bu_t)+\bs\e_t^{*}$ in (\ref{5}) to obtain
\begin{align}\label{nb1}
 \bs&\|\e_t^{*}\|^2+\frac{\nu}{2}\frac{d}{dt}\left(\bs\|\nabla\e^{*}\|^2\right)=\bs(\e_t^{*},\bu_t-P_h \bu_t)+\bs\nu\,a(\e^{*},\bu_t-P_h\bu_t)
+\frac{\nu}{2}\bs_t \|\nabla\e^{*}\|^2\nonumber\\
&-\bs b(\bu_H,\e^{*},P_h\e_t^{*})-\bs b(\e^{*},\bu_H,P_h\e_t^{*})+\bs b(\e_H,\e_H,P_h\e_t^{*})+\bs(p, \nabla\cdot P_h\e_t^{*}). 
\end{align}
Now, rewrite
\begin{align}\label{nb2}
 \bs b(\e_H,\e_H,P_h\e_t^{*})&=\frac{d}{dt}\left(\bs b(\e_H,\e_H,P_h\e^{*})\right)-\bs_t b(\e_H,\e_H,P_h\e^{*})-
\bs b(\e_{Ht},\e_H,P_h\e^{*})\nonumber\\
&-\bs b(\e_{H},\e_{Ht},P_h\e^{*}),
\end{align}
and
\begin{align}\label{nb2*}
 \bs(p, \nabla\cdot P_h\e_t^{*})&=\frac{d}{dt}\left(\bs (p-j_h p, \nabla\cdot P_h\e^{*})\right)-\bs_t (p-j_h p,\nabla\cdot P_h\e^{*})\nonumber\\
&-\bs (p_t-j_h p_t,\nabla\cdot P_h \e^{*}).
\end{align}
An application of Lemma \ref{nlt} with (\ref{3.4}) leads to 
\begin{align}\label{nb3}
 |\bs_t &b(\e_H,\e_H,P_h\e^{*})+\bs b(\e_{Ht},\e_H,P_h\e^{*})+\bs b(\e_H,\e_{Ht},P_h\e^{*})|\nonumber\\
&\leq C
\left(\bs_t\|\e_H\|^{1-\delta}\|\nabla\e_H\|^{1+\delta}+\bs\|\nabla\e_{Ht}\|\|\e_H\|^{1-\delta}\|\nabla\e_H\|^{\delta}\right)\|\nabla \e^{*}\|.
\end{align}
Apply (\ref{nb2})-(\ref{nb3}) along with the Cauchy-Schwarz's inequality, the Young's inequality, (\ref{3.4}) and Lemma \ref{nlt} in (\ref{nb1}). Then
integrate the resulting equation with respect to time from 0 to $t$ to arrive at
\begin{align}
 \int_0^t &\bs(\tau)\|\e_{\tau}^{*}(\tau)\|^2d\tau +\nu\,\bs\|\nabla\e^{*}\|^2\leq 
C(\nu)\bigg(h^2\int_0^t (\bs(\tau)\|\nabla\bu_t(\tau)\|^2+\bs_{\tau}(\tau)\|\nabla p(\tau)\|^2) d\tau\nonumber\\
&+h^2\int_0^t\frac{\bs^2(\tau)}{\bs_\tau(\tau)}\mathcal{K}^2_\tau(\tau)d\tau
+\int_0^t\left(e^{2\alpha \tau}\|\nabla\e^{*}(\tau)\|^2+\bs(\tau)\|\nabla\bu_H(\tau)\|\|\tilde\Delta\bu_H(\tau)\|
\|\nabla\e^{*}(\tau)\|^2\right)d\tau\nonumber\\
&+\bs (\|\e_H\|^{2(1-\delta)}\|\nabla\e_H\|^{2(1+\delta)}+h^2\|\nabla p\|^2)+\int_0^t\bs_\tau(\tau)\|\e_H(\tau)\|^{2(1-\delta)}
\|\nabla\e_H(\tau)\|^{2(1+\delta)}d\tau\nonumber\\
&+\int_0^t\frac{\bs^2(\tau)}{\bs_\tau(\tau)}
\|\nabla\e_{H\tau}(\tau)\|^2\|\e_H(\tau)\|^{2(1-\delta)}\|\nabla\e_H(\tau)\|^{2\delta}d\tau\bigg).\nonumber 
\end{align}
A use of $\frac{\bs(\tau)}{\bs_\tau(\tau)}\leq C \min\{1,\tau\}$ with $\bs_\tau(\tau)\leq Ce^{2\alpha \tau}$ 
and $\min\{1,\tau\}\leq \min\{1,t\}$ leads to
\begin{align}
 \int_0^t &\bs(\tau)\|\e_{\tau}^{*}(\tau)\|^2d\tau +\nu\,\bs\|\nabla\e^{*}\|^2\leq 
C(\nu)\bigg(h^2\int_0^t(\bs(\tau)\|\nabla\bu_t(\tau)\|^2+\bs_{\tau}(\tau)\|\nabla p(\tau)\|^2) d\tau\nonumber\\
&+h^2\min\{1,t\}\int_0^t\bs(\tau)\mathcal{K}^2_\tau(\tau)d\tau
+\int_0^te^{2\alpha \tau}\|\nabla\e^{*}(\tau)\|^2d\tau\nonumber\\
&+\int_0^t\bs(\tau)\|\nabla\bu_H(\tau)\|\|\tilde\Delta\bu_H(\tau)\|
\|\nabla\e^{*}(\tau)\|^2d\tau+\bs (\|\e_H\|^{2(1-\delta)}\|\nabla\e_H\|^{2(1+\delta)}+h^2\|\nabla p\|^2)	\nonumber\\
&+\int_0^t\bs_\tau(\tau)\|\e_H(\tau)\|^{2(1-\delta)}
\|\nabla\e_H(\tau)\|^{2(1+\delta)}d\tau\nonumber\\
&+\min\{1,t\}\int_0^t\bs(\tau)
\|\nabla\e_{H\tau}(\tau)\|^2\|\e_H(\tau)\|^{2(1-\delta)}\|\nabla\e_H(\tau)\|^{2\delta}d\tau\bigg).\nonumber 
\end{align}
From Lemmas \ref{L45}, \ref{T2}, \ref{nl5} and boundedness of $\|\nabla\bu_H\|$, it follows that
\begin{align}\label{nb5}
 \int_0^t &\bs(\tau)\|\e_{\tau}^{*}(\tau)\|^2d\tau +\nu\,\bs\|\nabla\e^{*}\|^2\leq C\bigg(h^2\bs+\bs 
\|\e_H\|_{L^{\infty}(\bL^2)}^{2(1-\delta)}\|\nabla\e_H\|_{L^{\infty}(\bL^2)}^{2(1+\delta)}\nonumber\\
&+\int_0^te^{2\alpha \tau}\|\nabla\e^{*}(\tau)\|^2d\tau+\|\e_H\|_{L^{\infty}(\bL^2)}^{-2\delta}\|\nabla\e_H\|_{L^{\infty}(\bL^2)}^{2(1+\delta)}
\int_0^t e^{2\alpha \tau}\|\e_H(\tau)\|^2d\tau\nonumber\\
&+\|\e_H\|_{L^{\infty}(\bL^2)}^{2-2\delta}\|\nabla\e_H\|_{L^{\infty}(\bL^2)}^{2\delta}\int_0^t\frac{\bs^2(\tau)}{\bs_\tau(\tau)}
\|\nabla\e_{H\tau}(\tau)\|^2d\tau\bigg).
 \end{align}
An application of Lemmas \ref{T2}, \ref{nl5} completes the proof.\hfill{$\Box$}

\begin{ldf}\label{nl6}
Under the hypotheses of Lemma \ref{T2}, the following estimate holds true:
 \begin{align}
   \bs^{-1}(t)\displaystyle{\int_0^t}e^{2\alpha \tau}\|\e^{*}(\tau)\|^2 d\tau
\leq K(t)\; (h^4+h^2H^4+H^6).\nonumber
\end{align} 
 \end{ldf}
\noindent
{\it Proof.}
Consider the linearized backward problem \cite{HS}. For a given $\e^{*}\in L^2(\bL^2)$, let
$(\bphi(t),\psi(t))\in \bJ_1\times L^2(\Omega)/\mathbb{R}$ be a weak solution of    
 \begin{align}\label{ne28}
  (\bv,\bphi_t)-\nu\,a(\bv,\bphi)-b(\bu,\bv,\bphi)-b(\bv,\bu,\bphi)+(\psi,\nabla\cdot\bv)
=(e^{2\alpha t}\e^{*},\bv),\,\,\forall \bv\in\bH_0^1(\Omega).
 \end{align}
with $\bphi(T)=0$, satisfying 
\begin{align}\label{ne29}
\displaystyle{\int_0^T}e^{-2\alpha t}(\|\bphi\|_2^2+\|\psi\|_1^2+\|\bphi_t\|^2)dt \leq 
C \displaystyle{\int_0^T}e^{2\alpha t}\|\e^{*}\|^2dt.
\end{align}
Rewrite (\ref{ne28}) as 
\begin{align}
 (\bv,\bphi_t)-\nu\,a(\bv,\bphi)-b(\bu_H,\bv,\bphi)-b(\bv,\bu_H,\bphi)-b(\e_H,\bv,\bphi)-b(\bv,\e_H,\bphi)
+(\psi,\nabla\cdot\bv)=(e^{\alpha t}\e^{*},\bv).\nonumber
\end{align}
Substitute $\bv=\e^{*}$, use (\ref{5}) with $\bphi_h=P_h\bphi$ and the discrete incompressibility condition to obtain
\begin{eqnarray}
e^{2\alpha t}\|\e^{*}\|^2
&=&\frac{d}{dt}(\e^{*},\bphi)-(\e^{*}_t,\bphi-P_h\bphi)-
\nu\,a(\e^{*},\bphi-P_h\bphi)-b(\bu_H,\e^{*},\bphi-P_h\bphi)\nonumber\\
&&-b(\e^{*},\bu_H,\bphi-P_h\bphi)-
b(\e_H,\e^{*},\bphi)-b(\e^{*},\e_H,\bphi)-b(\e_H,\e_H,P_h\bphi)\nonumber\\
&&-(p-j_h p,\nabla\cdot (P_h\bphi-\bphi))+(\psi-j_h\psi,\nabla\cdot\e^{*}).\nonumber
\end{eqnarray}
Using the fact
\begin{align}
 (\e^{*}_t,\bphi-P_h\bphi)&=\frac{d}{dt}(\e^{*}, \bphi-P_h\bphi)-(\e^{*}, \bphi_t-P_h\bphi_{t})\nonumber\\
&=\frac{d}{dt}(\e^{*},\bphi-P_h \bphi)-(\bu-P_h \bu,\phi_t),\nonumber
\end{align}
we now arrive at
\begin{eqnarray}
 e^{2\alpha t}\|\e^{*}\|^2
&=&\frac{d}{dt}(\e^{*},P_h\bphi)-(\bu-P_h\bu,\bphi_t)-\nu,a(\e^{*},\bphi-P_h\bphi)-b(\bu_H,\e^{*},\bphi-P_h\bphi)\nonumber\\
&&-b(\e^{*},\bu_H,\bphi-P_h\bphi)-b(\e_H,\e^{*},\bphi)-b(\e^{*},\e_H,\bphi)-b(\e_H,\e_H,P_h\bphi-\bphi)\nonumber\\
&&-b(\e_H,\e_H,\bphi)-(p-j_h p,\nabla\cdot (P_h\bphi-\bphi))+(\psi-j_h\psi,\nabla\cdot\e^{*}).\nonumber
\end{eqnarray}
Integrate with respect to time from $0$ to $T$ and use $\bphi(T)=0$ to find that
\begin{align}\label{ne34}
\displaystyle{\int_0^T} e^{2\alpha t}&\|\e^{*}(\tau)\|^2d\tau
=-(\e^{*}(0),P_h\bphi(0))-\displaystyle{\int_0^T}(\bu-P_h\bu,\bphi_\tau)d\tau\nonumber\\
&-\nu\displaystyle{\int_0^T}a(\e^{*},\bphi-P_h\bphi)d\tau-\displaystyle{\int_0^T}\left(b(\bu_H,\e^{*},\bphi-P_h\bphi)
+b(\e^{*},\bu_H,\bphi-P_h\bphi)\right)d\tau\nonumber\\
&-\displaystyle{\int_0^T}(\left(b(\e_H,\e^{*},\bphi)+b(\e^{*},\e_H,\bphi)+b(\e_H,\e_H,P_h\bphi-\bphi)
+b(\e_H,\e_H,\bphi)\right)d\tau\nonumber\\
&+\displaystyle{\int_0^T}\left(-(p-j_h p,\nabla\cdot (P_h\bphi-\bphi))
+(\psi-j_h\psi,\nabla\cdot\e^{*})\right) d\tau.
\end{align}
The first term in the right hand side of (\ref{ne34}) vanishes due to the orthogonality property of $P_h$. 
An application of (\ref{3.4}) with Cauchy-Schwarz's inequality and Young's inequality yields
\begin{align}\label{ne36}
 \int_0^T |(\bu-P_h\bu,\bphi_\tau)+a(\e^{*},\bphi-P_h\bphi)|d\tau &\leq C(\epsilon) 
 h^2\int_0^T e^{2\alpha \tau}\left(h^2\|\tilde\Delta \bu\|^2 +\|\nabla \e^{*}\|^2\right)d\tau\nonumber\\
 &+\epsilon \int_0^T e^{-2\alpha \tau}\|\bphi_\tau\|^2d\tau 
\end{align}
A use of (\ref{3.4}) with Lemma \ref{nlt} and boundedness of $\|\nabla\bu_H\|$  shows 
\begin{align}
 \int_0^T|b(\bu_H&,\e^{*},\bphi-P_h\bphi)+b(\e^{*},\bu_H,\bphi-P_h\bphi)|+|b(\e_H,\e^{*},\bphi)
 +b(\e^{*},\e_H,\bphi)|d\tau\nonumber\\
 &\leq C(\epsilon) h^2
\int_0^T e^{2\alpha \tau}(h^2+\|\e_H\|^2)\|\nabla\e^{*}\|^2d\tau
+\epsilon\int_0^Te^{-2\alpha \tau}\|\bphi\|_2^2d\tau.
\end{align}
Apply (\ref{3.4}) and \rm{\bf(B2)} to obtain
\begin{align}\label{ne39*}
 \int_0^T|(p-j_h p,&\nabla\cdot(P_h\bphi-\bphi))|+|(\psi-j_h\psi,\nabla\cdot\e^{*})|d\tau\leq 
C \int_0^T(h^2\|\nabla p\|\|\bphi\|_2+h\|\nabla\e^{*}\|\|\psi\|_1)d\tau\nonumber\\
&\leq C(\epsilon)\int_0^T e^{2\alpha \tau}\left(h^4\|\nabla p\|^2+ h^2\|\nabla\e^{*}\|^2\right)d\tau
+ \epsilon \int_0^T e^{-2\alpha \tau}\|\bphi\|_2^2d\tau. 
\end{align}
A use of (\ref{3.4}) with Lemma \ref{nlt} leads to
\begin{align}\label{ne38}
 \int_0^T (|b(\e_H,&\e_H,P_h\bphi-\bphi)|+|b(\e_H,\e_H,\bphi)|)d\tau
\leq C(\epsilon)\big(h^2\displaystyle{\int_0^T} e^{2\alpha \tau}\|\e_H\|^{2(1-\delta)}\|\nabla\e_H\|^{2(1+\delta)}d\tau\nonumber\\
&+\displaystyle{\int_0^T} e^{2\alpha \tau}\|\e_H\|^2\|\nabla\e_H\|^2d\tau\big)+\epsilon \int_0^T e^{-2\alpha \tau}\|\bphi\|_2^2d\tau.
\end{align}
Substitute (\ref{ne36})-(\ref{ne38}), regularity results (\ref{ne29}) and Lemma \ref{nl5} in (\ref{ne34}) 
with $\epsilon=1/4$ to obtain
\begin{align}\label{ne38*}
\displaystyle{\int_0^T} e^{2\alpha \tau}&\|\e^{*}(\tau)\|^2d\tau 
\leq C \bigg(h^4\int_0^T e^{2\alpha \tau}\mathcal{K}^2(\tau)d\tau+
\int_0^T e^{2\alpha \tau}\|\nabla\e^{*}(\tau)\|^2(h^2+\|\e_H(\tau)\|^2)d\tau\nonumber\\
&h^2\displaystyle{\int_0^T} e^{2\alpha \tau}\|\e_H(\tau)\|^{2(1-\delta)}\|\nabla\e_H(\tau)\|^{2(1+\delta)}d\tau+\int_0^T e^{2\alpha \tau}\|\e_H(\tau)\|^2\|\nabla\e_H(\tau)\|^2d\tau \bigg)\nonumber\\
&\leq C \bigg(h^4\int_0^T e^{2\alpha \tau}\mathcal{K}^2(\tau)d\tau+\left(h^2+\|\e_H\|^2_{L^{\infty}(\bL^2)}\right)
\int_0^T e^{2\alpha \tau}\|\nabla\e^{*}(\tau)\|^2d\tau\nonumber\\
&\qquad+\|\nabla\e_H\|^2_{L^{\infty}(\bL^2)}\int_0^T e^{2\alpha \tau}
\|\e_H(\tau)\|^2d\tau \bigg).
\end{align}
 A use of Lemmas \ref{L45}, \ref{T2} and \ref{nl5} in (\ref{ne38*}) concludes the proof. \hfill{$\Box$}
 
The following theorem provides estimates for $\e^{*}$. 
\begin{tdf}\label{T7}
Let the assumptions of Theorem \ref{T5} be satisfied. Further, let the discrete initial intial velocity 
$\bu_{0h}^{*}\in \bJ_h$ with $\bu_{0h}^{*}=P_h\bu_0$, where $\bu_0\in\bJ_1$. 
Then, there exists a positive constant $C$, independent of $h$, such that for $t\in (0,T]$ with $0<T<\infty$, 
the following estimates hold true:
\begin{align}
\| (\bu-\bu_h^{*})(t)\|\leq K(t)(h^2+ H^{3-\delta}\big),\,\,
\| \nabla(\bu-\bu_h^{*})(t)\|\leq K(t)(h+ H^{3-\delta}\big)\nonumber,
\end{align}
where $\delta>0$ is arbitrarily small and $K(t)=C e^{Ct}$. If, in addition, uniqueness condition (\ref{tn}) holds true, then $K(t)=K$, that is, 
estimates are bounded uniformly with respect to time.
\end{tdf}
\noindent
{\it Proof.} Since $\e^{*}=\bzeta+\brho$ and estimates of $\bzeta$ are known from  Lemma \ref{nl2}, it is enough to  derive estimates of $\brho$. A use of (\ref{5}) with (\ref{E616*}) and 
(\ref{tw9*}) leads to
\begin{align}\label{ne41}
(\brho_t,{\bphi}_h)+\nu\,a(\brho,\bphi_h) &+b(\bu_H,\brho,\bphi_h)+b(\brho,\bu_H,\bphi_h)
=-(\bzeta_{t},\bphi_h)\nonumber\\
&+b(\e_H,\e_H,\bphi_h)\,\,\,\,{\rm for~ all}~ \bphi_h\in \bJ_h.
\end{align}
Multiplying (\ref{ne41}) by $\bs(t)$, substitute $\bphi_h=\brho$ and use (\ref{tw11}) to arrive at 
\begin{align}
 \frac{1}{2}\frac{d}{dt}(\bs \|\brho\|^2)+\gamma\bs\|\nabla\brho\|^2\leq -\bs~(\bzeta_{t},\brho)+\bs~ 
b(\e_H,\e_H,\brho)+\frac{1}{2}\bs_t \|\brho\|^2,\nonumber 
\end{align}
where $\bs(t)=\rm {min}\{t,1\}~e^{2\alpha t}$.\\
Integrate with respect to time from $0$ to $t$ and obtain
\begin{align}\label{ne44}
 \bs(t)\|\brho(t)\|^2&+\gamma\int_0^t \bs(\tau)\|\nabla\brho(\tau)\|^2 d\tau\leq-\int_0^t\bs(\tau)
(\bzeta_{\tau}(\tau),\brho)d\tau\nonumber\\
&+\int_0^t\bs(\tau)\,b(\e_H(\tau),\e_H(\tau),\brho) d\tau
+\frac{1}{2}\int_0^t\bs_\tau(\tau) \|\brho(\tau)\|^2d\tau.
\end{align}
Using the Cauchy-Schwarz inequality and Lemma \ref{nl2}, the first term on the right hand 
side of (\ref{ne44}) can be bounded as 
 \begin{align}\label{ne45}
 & \int_0^t\bs(\tau)(\bzeta_\tau(\tau),\brho)d\tau\leq\int_0^t \bs(\tau)\|\bzeta_\tau(\tau)\|\|\brho\|d\tau\nonumber\\
&\leq \int_0^t \bs(\tau)h(h+\|\e_H(\tau)\|)(\mathcal{K}(\tau)\|\nabla\bu_{H\tau}(\tau)\|+\mathcal{K}_\tau (\tau))\|\brho\|d\tau. 
 \end{align}
A use of The Young's inequality with estimates of $\|\nabla\bu_{Ht}\|$, Lemmas \ref{L45} and \ref{T2} in (\ref{ne45}) leads to
\begin{align}\label{new1}
 \int_0^t\bs(\tau)&(\bzeta_\tau(\tau),\brho)d\tau\leq C(\lambda_1,\epsilon)h^2(h^2+\|\e_H\|_{L^{\infty}(\bL^2)}^2)
\bigg(\displaystyle{\sup_{0<t<T}}(\tau(t)\|\nabla\bu_{Ht}\|^2)\int_0^t
e^{2\alpha \tau}\mathcal{K}^2(\tau)d\tau\nonumber\\
&+\int_0^t \frac{\bs^2(\tau)}{\bs_\tau(\tau)}
\mathcal{K}^2_\tau(\tau)d\tau\bigg)+ \int_0^t \bs_\tau(\tau)\|\brho(\tau)\|^2d\tau+\epsilon \int_0^t
\bs(\tau)\|\nabla\brho(\tau)\|^2\;d\tau.
\end{align}
For the second term in the right hand side of (\ref{ne45}), split $\rho=\e^{*}-\bzeta$ and use Lemmas \ref{L45} and \ref{T2} to obtain
\begin{align}\label{ne47}
\int_0^t \bs_\tau(\tau)&\|\brho(\tau)\|^2d\tau\leq \int_0^t e^{2\alpha \tau}
(\|\e^{*}(\tau)\|^2+\|\bzeta(\tau)\|^2)d\tau\nonumber\\
&\leq C\bigg(\left(h^4+h^2\|\e_H\|^2_{L^{\infty}(\bL^2)}\right)\int_0^T e^{2\alpha \tau}\mathcal{K}^2(\tau)d\tau+
\|\nabla\e_H\|^2_{L^{\infty}(\bL^2)}\int_0^T e^{2\alpha \tau}\|\e_H(\tau)\|^2d\tau\bigg).
\end{align}
A use of Lemmas \ref{nlt} and \ref{T2} yields
\begin{align}\label{ne48}
 \int_0^t\bs(\tau) & b(\e_H(\tau),\e_H(\tau),\brho)d\tau\leq C \int_0^t \bs(\tau)\|\e_H(\tau)\|^{1-\delta}\|\nabla\e_H(\tau)\|^{1+\delta}
\|\nabla \brho\|d\tau \nonumber\\
&\leq C(\epsilon)\|\e_H\|_{L^{\infty}(\bL^2)}^{-2\delta}\|\nabla\e_H\|_{L^{\infty}(\bL^2)}^{2(1+\delta)}\int_0^t \bs(\tau)
\|\e_H(\tau)\|^2d\tau+\epsilon\,\int_0^t\bs(\tau)
\|\nabla\brho(\tau)\|^2\;d\tau.
\end{align}
An application of (\ref{ne45})-(\ref{ne48}) in (\ref{ne44}) leads to
\begin{align}\label{ne49}
 \|\brho(t)\|^2+\bs^{-1}(t)\int_0^t \bs(\tau)\|\nabla\brho(\tau)\|^2d\tau &\leq  K(t) \left(h^4+H^{6-2\delta}\right).
\end{align}
A use of triangle inequality with (\ref{ne49}) and Lemmas \ref{nl2}, \ref{nl5*} 
completes the proof of Theorem \ref{T7}. For establishing uniform estimates in Theorem \ref{T7}, 
use Lemma \ref{T2} in (\ref{new1})-(\ref{ne48}).  \hfill{$\Box$}

Next, we derive the error estimate for the two-grid approximation $p_h^{*}$ of the pressure $p$. Now, 
consider an equivalent form of (\ref{2a1}), that is, find $(\bu_h^{*}(t),p_h^{*}(t)) \in {\bf H}_h \times L_h$ 
such that $ \bu_h^{*}(0)= \bu_{0h} $ and for $t>0$
 \begin{eqnarray}\label{ng1}
 \left.
 \begin{array}{rcl}
&&(\bu_{ht}^{*},\bphi_h)+\nu a(\bu_h^{*}, \bphi_h) 
+b(\bu_h^{*},\bu_H,\bphi_h)+b(\bu_H,\bu_h^{*},\bphi_h)\\
&&\hspace{1.5cm}=(\f,\bphi_h)+b(\bu_H,\bu_H,\bphi_h)
+(p_h^{*}, \nabla \cdot {\bphi}_h)\;\;\; \forall {\bphi}_h \in {\bf H}_h,\\
&&(\nabla \cdot \bu_h^{*}, \chi_h) = 0 \;\;\; \forall \chi_h \in L_h.
\end{array}
 \right\}
 \end{eqnarray} 
To estimate $p-p_h^{*}$, use $j_h p$ and triangle inequality to obtain 
\begin{align}\label{ng2}
 \|p-p_h^{*}\|\leq \|p-j_h p\|+\|j_h p-p_h^{*}\|.
\end{align}
 From $({\bf B2})$, observe that
 \begin{align}\label{ng3}
\|j_h p-p_h^{*}\|_{L^2/N_h}&\leq C \displaystyle{\sup_{{\bphi}_h\in\bH_h\diagdown\{0\}}}{\left\{\frac{|(j_h p-p_h^{*},
 \nab\cdot{\bphi}_h)|}{\|\nab{\bphi}_h\|}\right\}}\nonumber\\
 &\leq C\bigg(\|j_h p-p\|+\displaystyle{\sup_{{\bphi}_h\in\bH_h\diagdown\{0\}}}{\left\{\frac{|( p-p_h^{*},
 \nab\cdot{\bphi}_h)|}{\|\nab{\bphi}_h\|}\right\}}\bigg).
 \end{align}
 The first term on the right hand side of (\ref{ng3}) can be estimated using $({\bf B1})$. 
 To estimate the second term on the right hand side of (\ref{ng3}), subtract (\ref{ng1}) from (\ref{4a1}) to obtain
\begin{align}\label{ng4}
(p-p_h^{*}&,\nabla \cdot {\bphi}_h)=(\e_t^{*},{\bphi}_h)+\nu\, a(\e^{*},{\bphi}_h)+b(\e^{*},\e_H,\bphi_h)\nonumber\\
&+b(\e_H,\e^{*},\bphi_h)-b(\e_H,\e_H,\bphi_h)\,\,\,\,\,\, \forall {\bphi}_h\in\bH_h.
\end{align}
A use of Lemma \ref{nlt} yields 
\begin{align}\label{ng5}
 |b(\e^{*}&,\e_H,\bphi_h)+b(\e_H,\e^{*},\bphi_h)-b(\e_H,\e_H,\bphi_h)|\nonumber\\
 &\leq C (\|\nabla\e_H\|\|\nabla\e^{*}\|+\|\e_H\|^{1-\delta}
\|\nabla\e_H\|^{1+\delta})\|\nabla\bphi_h\|.
\end{align}
Apply the Cauchy-Schwarz's inequality with (\ref{ng5}) to arrive at 
\begin{align}\label{ng6}
(p-p_h^{*},\nab\cdot{\bphi}_h)\le C\big(\|{\bf
e}^{*}_t\|_{-1;h}+ \|\nab{\bf e}^{*}\|+\|\nabla\e_H\|\|\nabla\e^{*}\|+\|\e_H\|^{1-\delta}
\|\nabla\e_H\|^{1+\delta}\big)\|\nab{\bphi}_h\|,
\end{align}
where 
\begin{align}
 \|{\bf e}^{*}_t\|_{-1;h}=\sup\left\{\frac{<\e^{*}_t,\bphi_h>}{\|\nabla\bphi_h\|}:\bphi_h\in\bH_h,\bphi_h\neq0\right\}.
\end{align}
Since the estimate of $\|\nab{\bf e}^{*}\|$ is known from Lemma \ref {nl5*}, we now derive 
estimate of  $\| {\bf e}^{*}_t\|_{-1;h}$. As $\bH_h\subset\bH_0^1$, we note that
\begin{align}\label{5n3}
 \|{\bf e}^{*}_t\|_{-1;h}&=\sup\left\{\frac{<\e^{*}_t,\bphi_h>}{\|\nabla\bphi_h\|}:\bphi_h\in\bH_h,\bphi_h\neq0\right\}\nonumber\\
 &\leq \sup\left\{\frac{<\e^{*}_t,\bphi>}{\|\nabla\bphi\|}:\bphi\in\bH_0^1,\bphi\neq0\right\}=\|\e_t^{*}\|_{-1}.
\end{align}
\begin{ldf}\label{nl8r}
 The error $\e^{*}=\bu-\bu_h^{*}$ satisfies for $0<t<T$
 \begin{align}
  \|\e_t^{*}\|_{-1}\leq K(t)\; (h+H^{3-\delta}).
 \end{align}
\end{ldf}
\noindent{\it Proof.} For any $\Psi\in \bH_0^1$, use orthogogal projection $P_h:\bL^2\rightarrow \bJ_h$ and 
(\ref{5}) with $\bphi=P_h\Psi$ to obtain
\begin{align}\label{5n1}
 (\e^{*}_t,\Psi)&=(\e_t^{*},\Psi-P_h\Psi)+(\e_t^{*},P_h\Psi)\nonumber\\
 &=(\e_t^{*},\Psi-P_h\Psi)-\nu\, a(\e^{*},P_h \Psi)-b(\e^{*},\bu_H,P_h\Psi)-b(\bu_H,\e^{*},P_h\Psi)\nonumber\\
 &+b(\e_H,\e_H,P_h\Psi)+(p,\nabla \cdot P_h\Psi).
\end{align}
Apply approximation property of $P_h$ to find that 
\begin{align}\label{5n2}
 (\e_t^{*},\Psi-P_h\Psi)=(\bu_t,\Psi-P_h\Psi)\leq C h\|\bu_t\|\|\nabla P_h\Psi\|.
\end{align}
A use of the Cauchy-Schwarz's inequality with (\ref{ng5}), (\ref{5n2}), the discrete incompressibility condition and (\ref{5n3}) in (\ref{5n1}) leads to
\begin{align}\label{5n4}
 \|\e_t^{*}\|_{-1}\leq C (h \|\bu_t\|+\|\nabla\e^{*}\|+\|\nabla\e_H\|\|\nabla\e^{*}\|+\|\e_H\|^{1-\delta}
\|\nabla\e_H\|^{1+\delta}+h\|\nabla p\|).
\end{align}
Using Lemmas \ref{L45}, \ref{T2} and \ref{nl5*}, 
we arrive at the desired result.\hfill{$\Box$}

Thus, we have proved the following theorem. 
\begin{tdf}\label{tp1n} 
Under the hypotheses of Lemma \ref{T2},  
there exists a positive constant $ C $ depending on $~ \nu,~ \gamma, ~M_0$, such that, for all $t>0$, it holds:
$$
\|(p-p_h^{*})(t)\|_{L^2/N_h}\leq K(t) (h+H^{3-\delta}),\nonumber
$$
where $\delta>0$ is arbitrarily small and $K(t)=C e^{Ct}$. Under the uniqueness condition (\ref{tn}), $K(t)=C$ and 
the estimate is uniform in time.
\end{tdf}
\subsection{Error estimates for Step 3}
This section is devoted to the derivation of semidiscrete error estimates in {\bf Step 3}. 
\begin{ldf}\label{nl8*}
Under the hypotheses of Lemma \ref{nl5}, the following estimate holds true:
\begin{align}
\tau(t)\|\e^{*}_t\|^2+\beta\;\bs^{-1}(t)\int_0^t\bs_1(\tau)\|\nabla\e_t^{*}(\tau)\|^2\leq \; K(t)\;\Big(h^2+H^{6-2\delta} \Big),\nonumber
\end{align}
where $\bs_1(t)=\tau^2(t)e^{2\alpha t}$ and $\beta= \gamma-\alpha\lambda_1^{-1}>0$.
\end{ldf}
\noindent{\it Proof.} 
Differentiate (\ref{5}) with respect to time and substitute $\bphi_h=P_h\e_t^{*}$ in the resulting equation and use 
discrete incompressibility condition to arrive at
\begin{align}\label{nd2}
\frac{1}{2}\frac{d}{dt}&\|\e_t^{*}\|^2+\gamma\|\nabla\e_t^{*}\|^2\leq (e_{tt}^{*}, \bu_t-P_h\bu_t)+
\nu\,a(\e_t^{*},\bu_t-P_h\bu_t)-b(\e_t^{*},\bu_{H},\bu_t-P_h\bu_t)\nonumber\\
&-b(\bu_{H},\e_t^{*},\bu_t-P_h\bu_t)-b(\e^{*},\bu_{Ht},P_h\e_t^{*})
-b(\bu_{Ht},\e^{*},P_h\e_t^{*})\nonumber\\
&+b(\e_{Ht},\e_H,P_h\e_t^{*})+b(\e_{Ht},\e_H,P_h\e_t^{*})+(p_t-j_h p_t,\nabla\cdot P_h\e_t^{*}).
\end{align}
A use of (\ref{3.4})with Lemma \ref{nlt} and Theorem \ref{T2} yields
\begin{eqnarray}\label{nd6}
 |b(\e_{Ht},\e_H,P_h\e_t^{*}) &+b(\e_{H},\e_{Ht},P_h\e_t^{*})|\leq C \|\e_H\|^{1-\delta}\|\nabla\e_H\|^{\delta}\|\nabla\e_{Ht}\|
\|\nabla P_h\e_t^{*}\|\nonumber\\
&\leq C(\epsilon)\; \|\e_H\|^{2(1-\delta)}\;\|\nabla\e_H\|^{2\delta} \; \|\nabla\e_{Ht}\|^2+\epsilon\|\nabla \e_t^{*}\|^2.
\end{eqnarray}
Apply (\ref{ne22})-(\ref{ne23*}) and (\ref{nd6}) along with Cauchy-Schwarz's inequality in (\ref{nd2}) to arrive at 
\begin{eqnarray}\label{nd7}
\frac{1}{2}\frac{d}{dt}\|\e_t^{*}\|^2+\gamma\|\nabla\e_t^{*}\|^2&\leq &\frac{1}{2}\frac{d}{dt}\|\bu_t-P_h \bu_t\|^2
+ C(\nu,\gamma)\bigg(h^2\mathcal{K}_t^2+ h^2\|\nabla\bu_H\|\|\tilde\Delta\bu_H\|\|\tilde\Delta\bu_t\|^2\nonumber\\
&&+K(t)\left((h^2+H^{6-2\delta})\|\nabla \bu_{Ht}\|^2+H^{4-2\delta}\|\nabla\e_{Ht}\|^2\right)\bigg).
\end{eqnarray}
Use similar analysis to (\ref{nd7}) as applied to (\ref{ne25}) to arrive at (\ref{ne25*}) and Lemmas \ref{L45}, \ref{T2}, 
\ref{nl5*} to conclude the proof. \hfill{$\Box$} 
\begin{ldf}\label{nl8}
Under the hypotheses of Lemma \ref{nl5}, the following estimate holds true:
\begin{align}
 \bs^{-1}(t)\int_0^t e^{2\alpha \tau}\|\nabla \e_h(\tau)\|^2d\tau\leq \; K(t) (h^2+H^{10-4\delta}).\nonumber
\end{align}
\end{ldf}
\noindent
{\it Proof.} Consider (\ref{7}) with $\bphi_h=P_h \e_h=(P_h \bu-\bu)+\e_h$. Then, use (\ref{tw11}) to arrive at
 \begin{align}\label{na1}
 \frac{1}{2}\frac{d}{dt}\|\e_h\|^2+\gamma\,\|\nabla\e_h\|^2&\leq (\e_{ht},\bu-P_h \bu)
 +\nu\,a(\e_h,\bu-P_h\bu)+b(\bu_H,\e_h,\bu-P_h\bu)\nonumber\\
&\qquad+b(\e_h,\bu_H,\bu-P_h\bu)-b(\e_H,\e^{*},P_h \e_h)
\nonumber\\
&\qquad-b(\e^{*},\e_H,P_h\e_h)
+b(\e^{*},\e^{*},P_h\e_h)+(p,\nabla\cdot P_h\e_h).
 \end{align} 
The first four terms in the right hand side of (\ref{na1}) can be bounded using (\ref{ne22})-(\ref{ne24*}) 
(with $\e^{*}$ replaced by $\e_h$). Now, from (\ref{3.4}) and Lemma \ref{nlt}, we arrive at
\begin{align}\label{na2}
 |b(\e_H,\e^{*},P_h \e_h)&-b(\e^{*},\e_H,P_h\e_h)+b(\e^{*},\e^{*},P_h\e_h)|\nonumber\\
&\leq C \left(\|\e_H\|^{1-\delta}\|\nabla\e_H\|^{\delta}\|\nabla\e^{*}\|+\|\nabla\e^{*}\|^2\right)
\|\nabla\e_h\|.
\end{align}
A use of (\ref{ne22})-(\ref{ne23*}) and (\ref{na2}) leads to
\begin{align}\label{na3}
 \frac{1}{2}\frac{d}{dt}&\|\e_h\|^2+\gamma\|\nabla \e_h\|^2
\leq  \frac{1}{2}\frac{d}{dt}\|\bu-P_h\bu\|^2+C(\nu)h\mathcal{K}\|\nabla\e_h\|\nonumber\\
&+C(\nu)(\|\e_H\|^{1-\delta}\|\nabla\e_H\|^{\delta}\|\nabla\e^{*}\|+\|\nabla\e^{*}\|^2)\|\nabla\e_h\|.
\end{align}
The proof can be concluded by using the similar set of arguments now to (\ref{na3}) as applied to (\ref{ne25}) leading to 
(\ref{ne25*}) and Lemmas \ref{T2}, \ref{nl5} and \ref{nl5*}. This completes the proof of Lemma \ref{nl8}.\hfill{$\Box$}

Below, we state a lemma that provides $L^{\infty}(\bH^1)$-norm estimate for $\e_{h}$. The proof is obtained in the same lines as 
the proof of Lemma \ref{nl5*}, starting with $\bphi_h=P_h\e_{ht}$ in ({\ref{7}), using Lemmas \ref{nl8*} and \ref{nl8} 
and hence, is skipped.
\begin{ldf}\label{nl5a}
 Under the hypotheses of Lemma \ref{nl5}, the following estimate holds true:
\begin{align}
\bs^{-1}(t)\int_0^t\bs(\tau)\|\e_{h\tau}(\tau)\|^2d\tau+ \|\nabla \e_h(t)\| \leq \;K(t)\, (h^2+H^{8-2\delta}),\nonumber
\end{align}
where $\delta>0$ is arbitrarily small. 
\hfill{$\Box$}
\end{ldf}
\begin{ldf}\label{nl7}
Under the hypotheses of Theorem \ref{T2}, the following is satisfied:
\begin{align}
 \bs^{-1}(t)\int_0^te^{2\alpha \tau}\|\e_h(\tau)\|^2d\tau\leq \;K(t)(h^4+h^2H^4).\nonumber
\end{align}
\end{ldf}
\noindent
{\it Proof.} Proceeding in a similar way as in the proof of Lemma \ref{nl6}, we arrive at
\begin{eqnarray}\label{ne54}
\|\e_h\|^2
&=&\frac{d}{dt}(\e_h,P_h\bphi)-(\bu-P_h\bu,\bphi_t)-a(\e_h,\bphi-P_h\bphi)-b(\bu_H,\e_h,\bphi-P_h\bphi)\nonumber\\
&&-b(\e_h,\bu_H,\bphi-P_h\bphi)-b(\e_H,\e_h,\bphi)-b(\e_h,\e_H,\bphi)+b(\e_H,\e^{*},P_h\bphi-\bphi)\nonumber\\
&&+b(\e_H,\e^{*},\bphi)+b(\e^{*},\e_H,P_h\bphi-\bphi)+b(\e^{*},\e_H,\bphi)-b(\e^{*},\e^{*},P_h\bphi-\bphi)\nonumber\\
&&-b(\e^{*},\e^{*},\phi)-(p-j_h p,\nabla\cdot (P_h\bphi-\bphi))+(\psi-j_h\psi,\nabla\cdot\e_h).
\end{eqnarray} 
Use (\ref{3.4}), Lemmas \ref{nlt} and \ref{T2} to bound
\begin{align}\label{ne54*}
 |b(\e_H&,\e^{*},P_h\bphi-\bphi)|+|b(\e^{*},\e_H,P_h\bphi-\bphi)|+|b(\e_H,\e^{*},\bphi)|+|b(\e_H,\e^{*},\bphi)|\nonumber\\
 &\leq C (h\|\e_H\|^{1-\delta}\|\nabla\e_H\|^{\delta}+\|\e_H\|)
\|\nabla\e^{*}\|\|\bphi\|_2
 \leq K(t)( h H^{2-\delta}+H^2)\|\nabla\e^{*}\|\|\bphi\|_2.
\end{align}
An application of (\ref{3.4}), Lemmas \ref{nlt} and \ref{nl5*} yields
\begin{align}\label{ne55*}
| b(\e^{*},\e^{*},P_h\bphi-\bphi)|+|b(\e^{*},\e^{*},\phi)|&\leq C (h\|\nabla \e^{*}\|^2+\|\e^{*}\|\|\nabla\e^{*}\|)\|\bphi\|_2\nonumber\\
 &\leq K(t)(h^2+H^{3-\delta})\|\nabla\e^{*}\|\|\bphi\|_2.
\end{align}
Multiply (\ref{ne54}) by $e^{2\alpha t}$ and integrate with respect to time from $0$ to $t$. Then, apply 
(\ref{ne36})-(\ref{ne39*}) with $\e^{*}$ replaced by $\e_h$ and (\ref{ne54*})-(\ref{ne55*}) to obtain
\begin{eqnarray}
 \label{new3}
\displaystyle{\int_0^t} e^{2\alpha \tau}\|\e_h(\tau)\|^2d\tau 
&\leq &C \big(h^4\int_0^t e^{2\alpha \tau}\mathcal{K}^2(\tau)d\tau+
\int_0^t e^{2\alpha \tau}\|\nabla\e_h(\tau)\|^2(h^2+\|\e_H(\tau)\|^2)d\tau\big)\nonumber\\
&&+K(t)(h^2+H^{4-2\delta})e^{Ct}\int_0^t e^{2\alpha \tau}\|\nabla\e^{*}(\tau)\|^2d\tau\nonumber\\
&\leq& C\big(h^4\int_0^t e^{2\alpha \tau}\mathcal{K}^2(\tau)d\tau+
(h^2+\|\e_H\|_{L^{\infty}(L^2)}^2)\int_0^t e^{2\alpha \tau}\|\nabla\e_h(\tau)\|^2d\tau\big)\nonumber\\
&&+K(t)(h^2+H^{4-2\delta})\int_0^t e^{2\alpha \tau}\|\nabla\e^{*}(\tau)\|^2d\tau.
\end{eqnarray}
A use of Lemmas \ref{L45}, \ref{T2}, \ref{nl5} and \ref{nl8} completes the rest part of the proof.\hfill{$\Box$}

Proceeding in a similar way as in Lemma \ref{nl8r}, we arrive at the following estimate.
\begin{ldf}\label{n1}
  The error $\e_h=\bu-\bu_h$ satisfies for $0<t<T$
 \begin{align}
  \|\e_{ht}\|_{-1}\leq\; K(t) (h+H^{4-\delta}).\nonumber
 \end{align}
\end{ldf}

\noindent
For the pressure error estimates corresponding to the correction in {\bf Step 3} of two-grid algorithm, 
consider the equivalent form of (\ref{3a1}):
seek $(\bu_h(t),p_h(t)) \in {\bf H}_h \times  L_h$ such that $ \bu_h(0)= \bu_{0h} $ and for $t>0$,
 \begin{eqnarray}\label{nh11}
 \left.
 \begin{array}{rcl}
&&(\bu_{ht},\bphi_h)+\nu a(\bu_h, \bphi_h) +b(\bu_h,\bu_H,\bphi_h)
+b(\bu_H,\bu_h,\bphi_h)=(\f,\bphi_h)\\
&&\hspace{1.5cm}+b(\bu_H,\bu_h^{*},\bphi_h)+b(\bu_h^{*},\bu_H-\bu_h^{*},\bphi_h)+(p_h,\nabla\cdot\bphi_h)\,\,\, \forall {\bphi}_h \in {\bf H}_h,\\
 &&(\nabla \cdot \bu_h, \chi_h) = 0 \;\;\; \forall \chi_h \in L_h.
\end{array}
 \right\}
 \end{eqnarray}

For pressure equation, subtract (\ref{nh11}) from (\ref{4a1}) to obtain
\begin{align}\label{nh12}
(p-p_h,\nab\cdot{\bphi}_h)&=(\e_{ht},{\bphi}_h)+\nu a(\e_h,{\bphi}_h)+b(\e_h,\bu_H,\bphi_h)+b(\bu_H,\e_h,\bphi_h)\nonumber\\
&-b(\e_H,\e^{*}, \bphi_h)-b(\e^{*},\e_H, \bphi_h)
-b(\e^{*},\e^{*},\bphi_h)\,\,\,\, \forall {\bphi}_h \in \bH_h.
\end{align}

Armed with these estimates, next we derive proof of main Theorem \ref{T5}.\\
\noindent
{\it Proof of Theorem \ref{T5}.}
Multiply (\ref{B18*}) by $\bs(t)$, substitute $\bphi_h=\Th$
 and integrate the resulting equation from $0$ to $t$ to obtain 
\begin{align}\label{ne50}
 \bs(t)\|&\Th(t)\|^2+\gamma\int_0^t \bs(\tau)\|\nabla\Th(\tau)\|^2 d\tau\leq -\int_0^t\bs(\tau)~(\bzeta_{\tau}(\tau),\Th)d\tau
+\int_0^t\bs_\tau(\tau)\|\Th(\tau)\|^2d\tau\nonumber\\
&+\int_0^t \bs(\tau)\left(-b(\e_H(\tau),\e^{*}(\tau),\Th)
-b(\e^{*}(\tau),\e_H(\tau),\Th)
+b(\e^{*}(\tau),\e^{*}(\tau),\Th)\right)d\tau.
\end{align}
The first term on the right hand side of (\ref{ne50}) can be tackled as in (\ref{ne45}). Write $\Th=\e_h-\bzeta$ and use Lemmas 
\ref{L45}, \ref{nl2} to obtain
\begin{align}\label{ne47*}
\int_0^t \bs_\tau(\tau)&\|\Th(\tau)\|^2d\tau\leq \int_0^t e^{2\alpha \tau}
(\|\e_h(\tau)\|^2+\|\bzeta(\tau)\|^2)d\tau \leq K(t) (h^4+h^2H^{4-2\delta})\bs.
\end{align}
Use Young's inequality and Lemmas \ref{nlt}, \ref{T2}, \ref{nl5} to bound 
\begin{align}\label{ne51}
 \big|\int_0^t &\bs(\tau)(b(\e_H(\tau),\e^{*}(\tau),\Th)+b(\e^{*}(\tau),\e_H(\tau),\Th))d\tau\big|\nonumber\\
&\leq C(\epsilon)\|\e_H\|_{L^{\infty}(\bL^2)}^{2-2\delta}\|\nabla\e_H\|_{L^{\infty}(\bL^2)}^{2\delta}
\int_0^t\bs(\tau)\|\nabla\e^{*}(\tau)\|^2d\tau
+\epsilon\int_0^t\bs(\tau)\|\nabla\Th(\tau)\|^2d\tau\nonumber\\
&\leq K(t)(h^2H^{4-2\delta}+H^{10-4\delta})\bs+\epsilon\int_0^t\bs(\tau)\|\nabla\Th(\tau)\|^2d\tau.
\end{align}
An application of Lemmas \ref{nlt}, \ref{nl5} and \ref{nl5*} leads to 
\begin{align}\label{ne52}
\int_0^t &\bs(\tau)b(\e^{*}(\tau),\e^{*}(\tau),\Th)d\tau\leq C \int_0^t \bs(\tau)\|\nabla\e^{*}(\tau)\|^2\|\nabla\Th\|d\tau\nonumber\\
&\leq C\|\nabla\e^{*}\|_{L^{\infty}(\bL^2)}^2\int_0^t \bs(\tau)\|\nabla\e^{*}(\tau)\|^2d\tau 
+\epsilon\int_0^t\bs(\tau)\|\nabla\Th\|^2d\tau\nonumber\\
& \leq K(t) (h^4+H^{12-4\delta})\bs+\epsilon\int_0^t\bs(\tau)\|\nabla\Th\|^2d\tau.
\end{align}
Apply Lemma \ref{T2} to obtain
\begin{align}\label{ne53}
 \|\Th(t)\|^2+\bs^{-1}(t)\int_0^t\bs(\tau)\|\nabla\Th(\tau)\|^2d\tau\leq \;K(t) (h^4+h^2H^{4-2\delta}+H^{10-2\delta}). 
\end{align}
A use of Lemmas \ref{nl2}, \ref{nl5a} with (\ref{ne53}) completes the proof of Theorem 
\ref{T5}.  


The uniform estimates in Theorem \ref{T5} can be achieved by using Lemma \ref{T2} under uniqueness condition.

For the pressure estimate (\ref{error-pressure}), 
a use of boundedness of $\|\nabla\bu_H\|$ and Lemmas \ref{nlt}, \ref{T2}, \ref{nl5*}, \ref{nl5a} leads to
\begin{align}\label{nh13}
 |b(\e_h,\bu_H,\bphi_h)&+b(\bu_H,\e_h,\bphi_h)-b(\e_H,\e^{*}, \bphi_h)-b(\e^{*},\e_H, \bphi_h)-b(\e^{*},\e^{*},\bphi_h)|\nonumber\\
&\leq C(\|\nabla\e_h\|\|\nabla\bu_H\|+\|\nabla\e_H\|\|\nabla\e^{*}\|+\|\nabla\e^{*}\|^2)\|\nabla\bphi_h\|.
\end{align}
A use of Cauchy-Schwarz's inequality and (\ref{nh13}) in (\ref{nh12}) leads to 
\begin{align}
(p-p_h,\nab\cdot{\bphi}_h)\leq C\left(\|\e_{ht}\|_{-1}+\|\nabla \e_{h}\|+\|\nabla\e_h\|\|\nabla\bu_H\|
+\|\nabla\e_H\|\|\nabla\e^{*}\|+\|\nabla\e^{*}\|^2\right) \|\nabla\bphi_h\|.\nonumber
\end{align}
A use of Lemmas \ref{T2}, \ref{nl5a}, \ref{n1} completes the proof of the pressure estimate (\ref{error-pressure}) and this concludes the rest of the proof of Theorem \ref{T5}.\hfill{$\Box$} 
\section{Backward Euler Method}
\setcounter{equation}{0}
For a complete discretization, we apply a backward Euler method for the time discretization. 
Let $\{t_n\}_{n=0}^{N}$ be a uniform partition of the time interval $[0,T]$ and $t_n=nk$, 
with time step $k>0$. For a sequence  $\{\bphi^n\}_{n\geq 0} \in \bJ_h$ defined on $[0, T]$, set $\bphi^n=\bphi(t_n)$, 
$\bar \partial _t \bphi^n=\left(\bphi^n-\bphi^{n-1}\right)/k$. \\
\noindent
The backward Euler method applied to (\ref{1a1})-(\ref{3a1}) is stated in terms of the following algorithm: \\
{\it {\bf Step 1.} Solve the nonlinear system (\ref{1.1}) on $\mathcal{T}_H$: find $\bU_H^n\in \bJ_H$, such 
that for all $\bphi_H\in \bJ_H$ for $\bU_H^{0}=P_H\bu_0$ and $t>0$
\begin{eqnarray}\label{1bd1}
(\bar\partial_t \bU_H^n,\bphi_H)+\nu\, a(\bU_H^n, \bphi_H) +b(\bU_H^n,\bU_H^n,\bphi_H)=({\bf f}^n,\bphi_H).
 \end{eqnarray}
{\bf Step 2.} Update on $\mathcal{T}_h$ with one Newton iteration: find $\bU^{n}\in \bJ_h$,
such that for all $\bphi_h\in \bJ_h$ for $\bU^{0}=P_h\bu_{0}$ and $t>0$
 \begin{align}\label{2bd1}
(\bar\partial_t\bU^n,\bphi_h)&+\nu\, a(\bU^n, \bphi_h) 
+b(\bU^n,\bU_H^n,\bphi_h)\nonumber\\
&+b(\bU^n_H,\bU^n,\bphi_h)=({\bf f}^n,\bphi_h)+b(\bU_H^n,\bU^n_H,\bphi_h).
\end{align}
{\bf Step 3.} Correct on $\mathcal{T}_h$: find $\bU^n_{h}\in\ \bJ_h$ such that, for all $\bphi_h\in \bJ_h$ 
for $\bU_h^0=P_h\bu_{0}$ and $t>0$
\begin{align}\label{3bd1}
(\bar\partial_t\bU_h^n,\bphi_h)+\nu\, a( &\bU_h^n, \bphi_h) +b(\bU_h^n,\bU_H^n,\bphi_h)
+b(\bU_H^n,\bU_h^n,\bphi_h)\nonumber\\
&=({\bf f}^n,\bphi_h)+b(\bU_H^n,\bU^{n},\bphi_h)+b(\bU^{n},\bU_H^n-\bU^{n},\bphi_h).
\end{align}
}
The results in Lemmas \ref{1l1}-\ref{bdlm2} will play an important role in the derivation of error estimates in this section.
\noindent
 \begin{ldf}\label{1l1}
Let $\bu_h^{*}$ be the solution of (\ref{2a1}) on some interval $[0,T),\,\,0<T<\infty$ satisfying 
$\bu_{0h}^{*}=P_h\bu_0$. 
Then, there exists a positive constant $C=C(\gamma,\nu,\alpha,\lambda_1, M_0)$, such that for\, 
$\displaystyle{0\leq\alpha< \frac{\gamma\lambda_1}{2}}$ for all $t>0$, the following holds true:
 \begin{align}
  \|\bu_h^{*}(t)\|^2+\|\nabla\bu_h^{*}(t)\|^2+e^{-2\alpha t}\displaystyle{\int_0^t}e^{2\alpha s}
  (\|\nabla\bu_h^{*}(s)\|^2+\|\tilde\Delta\bu_h^{*}(s)\|^2)ds\leq C.\nonumber
 \end{align}
\end{ldf}
\noindent
{\it Proof.} Multiply (\ref{2a1}) by $e^{\alpha t}$ for some $\alpha >0$ and set 
$\hat\bu_h^{*}=e^{\alpha t}\bu_h^{*}$. Substitute $\bphi_h=\hat\bu_h^{*}$ and use (\ref{2.1*}) 
($\|\hat\bu_h^{*}\|^2\leq \lambda_1^{-1} \|\nabla\hat\bu_h^{*}\|^2$) and (\ref{tw11}) to obtain
\begin{align}\label{e2}
\frac{1}{2} \frac{d}{dt}\|\hat\bu_h^{*}\|^2+\left(\gamma-\frac{\alpha}{\lambda_1}\right)\|\nabla\hat\bu_h^{*}\|^2
\leq (\hat\f,\hat\bu_h^{*})+e^{-\alpha t}b(\hat\bu_H,\hat\bu_H,\hat\bu_h^{*}).
\end{align}
An application of Cauchy-Schwarz's inequality and Young's inequality leads to
\begin{align}\label{e3}
 |(\hat \f,\hat\bu_h^{*})|\leq C(\lambda_1,\epsilon)\|\hat\f\|^2+\epsilon \|\nabla\hat\bu_h^{*}\|^2.
\end{align}
A use of Lemma \ref{nlt}, $\|\nabla\bu_H\|\leq C$ and Young's inequality yields
\begin{align}\label{e4}
 |e^{-\alpha t}b(\hat\bu_H,\hat\bu_H,\hat\bu_h^{*})|&\leq C e^{-\alpha t}\|\nabla\hat\bu_H\|^2\|\nabla\hat\bu_h^{*}\|\nonumber\\
&\leq C(\epsilon)\|\nabla\hat\bu_H\|^2+\epsilon \|\nabla\hat\bu_h^{*}\|^2. 
\end{align}
Apply (\ref{e3})-(\ref{e4}) in (\ref{e2}) with $\epsilon=\frac{\gamma}{2}$ and integrate 
the resulting equation with respect to time to obtain
\begin{align}
 \|\hat\bu_h^{*}(t)\|^2+\left(\gamma-\frac{2\alpha}{\lambda_1}\right)\int_0^t\|\nabla\hat\bu_h^{*}(s)\|^2ds\leq \|\bu_0\|^2+
C\displaystyle{\int_0^t}\left(\|\nabla\hat\bu_H(s)\|^2+\|\hat\f(s)\|^2\right)ds.\nonumber
\end{align}
Multiply above equation by $e^{-2\alpha t}$, use assumption {\bf (A2)} and the fact that
 $e^{-2\alpha t}\displaystyle{\int_0^t}e^{2\alpha s}ds =\frac{1}{2\alpha}(1-e^{-2\alpha t})$
to arrive at
\begin{align}\label{e1}
 \|\bu_h^{*}(t)\|^2+e^{-2\alpha t}\int_0^te^{2\alpha s}\|\nabla\bu_h^{*}(s)\|^2ds\leq C.
 \end{align}
\noindent
Next, multiply (\ref{2a1}) by $e^{\alpha t}$ and rewrite it as
\begin{align}\label{e7}
 (\hat\bu_{ht}^{*},\bphi_h)&-\nu a(\tilde\Delta_h\hat\bu_h^{*},\bphi_h)=\alpha (\hat\bu_h^{*},\bphi_h)-
e^{-\alpha t}b(\hat\bu_H,\hat\bu_h^{*},\bphi_h)\nonumber\\
&-e^{-\alpha t} \left(b(\hat\bu_h^{*},\hat\bu_H,\bphi_h)
-b(\hat\bu_H,\hat\bu_H,\bphi_h)\right)+(\hat\f,\bphi_h).
\end{align}
Substitute $\bphi_h=-\tilde\Delta_h\hat\bu_h^{*}$ in (\ref{e7}), 
note the fact that $-(\hat\bu_{ht}^{*},\tilde\Delta\hat\bu_h^{*})=
\frac{1}{2}\frac{d}{dt}\|\nabla\hat\bu_h^{*}\|^2$ and integrate the resulting equation with respect to time to obtain
\begin{align}\label{e8}
\|\nabla\hat\bu_h^{*}(t)\|^2&+2\nu\displaystyle{\int_0^t}\|\tilde\Delta\hat\bu_h^{*}(s)\|^2 ds=\|\nabla\bu_{0h}^{*}\|^2-
2\alpha \displaystyle{\int_0^t}
(\hat\bu_h^{*},\tilde\Delta_h\hat\bu_h^{*})ds+
2\displaystyle{\int_0^t}e^{-\alpha s}b(\hat\bu_H,\hat\bu_h^{*},\tilde\Delta_h\hat\bu_h^{*})ds\nonumber\\
&+2\displaystyle{\int_0^t} e^{-\alpha s} \left(b(\hat\bu_h^{*},\hat\bu_H,\tilde\Delta_h\hat\bu_h^{*})
-b(\hat\bu_H,\hat\bu_H,\tilde\Delta_h\hat\bu_h^{*})\right)ds-2\displaystyle{\int_0^t}(\hat\f,\tilde\Delta_h\hat\bu_h^{*})ds.
\end{align}
An application of Lemmas \ref{nlt}, \ref{1l1*} with Young's inequality yields
\begin{align}\label{e9}
 2\displaystyle{\int_0^t}e^{-\alpha s}&\left(|b(\hat\bu_H,\hat\bu_h^{*},\tilde\Delta_h\hat\bu_h^{*})|
+|b(\hat\bu_h^{*},\hat\bu_H,\tilde\Delta_h\hat\bu_h^{*})|+|b(\hat\bu_H,\hat\bu_H,\tilde\Delta_h\hat\bu_h^{*})|\right)ds\nonumber\\
&\leq  C(\epsilon)\displaystyle{\int_0^t}e^{-2\alpha s}(\|\tilde\Delta_H \hat\bu_H(s)\|^2 +\|\nabla\hat\bu_h^{*}(s)\|^2)ds
+\epsilon\displaystyle{\int_0^t}\|\tilde\Delta_h\hat\bu_h^{*}(s)\|^2ds.
\end{align}
A use of (\ref{e9}) along with Cauchy-Schwarz's inequality leads to
\begin{align}
 \|\nabla\hat\bu_h^{*}(t)\|^2+\nu\displaystyle{\int_0^t}\|\tilde\Delta_h\hat\bu_h^{*}(s)\|^2 ds\leq \|\nabla\bu_{0h}^{*}\|^2
+C\displaystyle{\int_0^t}(\|\hat\f(s)\|^2+\|\nabla\hat\bu_h^{*}(s)\|^2+\|\tilde\Delta_H \hat\bu_H(s)\|^2)ds.\nonumber
\end{align}
An application of (\ref{e1}), assumption {\bf (A2)} and Lemma \ref{1l1*} completes the proof.\hfill{$\Box$}
\begin{ldf}\label{1l2}
Under the assumption of Lemma \ref{1l1}, the following holds true:
\begin{align}
 e^{-2\alpha t}\displaystyle{\int_0^t}e^{2\alpha s}(\|\bu_{ht}^{*}(s)\|^2+\|\bu_{htt}^{*}(s)\|^2_{-1})ds\leq C.\nonumber
\end{align}
\end{ldf}
\noindent
{\it Proof.} Substitute $\bphi_h=e^{2\alpha t}\bu_{ht}^{*}$ in (\ref{2a1}) and write it as 
\begin{align}
 e^{2\alpha t}\|\bu_{ht}^{*}\|^2&=\nu\,e^{2\alpha t}(\tilde\Delta_h\bu_h^{*},\bu_{ht}^{*})-
e^{2\alpha t}\left(b(\bu_{H},\bu_{h}^{*},\bu_{ht}^{*})+b(\bu_h^{*},\bu_H,\bu_{ht}^{*})\right)\nonumber\\
&+e^{2\alpha t}b(\bu_H,\bu_H,\bu_{ht}^{*})+e^{2\alpha t}(\f,\bu_{ht}^{*}).
\end{align}
Apply Lemmas \ref{nlt}, \ref{1l1*}, \ref{1l1} with Cauchy-Schwarz's inequality and Young's inequality to obtain
\begin{align}\label{e19}
  e^{2\alpha t}\|\bu_{ht}^{*}\|^2\leq Ce^{2\alpha t}(\|\tilde\Delta_H\bu_H\|^2+\|\tilde\Delta\bu_h^{*}\|^2+\|\f\|^2).
\end{align}
Integrate (\ref{e19}) with respect to time and use Lemmas \ref{1l1*} and \ref{1l1}, assumption {(\bf A2)} to arrive at 
\begin{align}\label{e16}
e^{-2\alpha t}\displaystyle{\int_0^t}e^{2\alpha s}\|\bu_{ht}^{*}(s)\|^2ds\leq C.
\end{align}
Next, differentiate (\ref{2a1}) with respect to time and obtain
\begin{align} \label{e11*}
 (\bu_{htt}^{*},\bphi_h)&+\nu\,a(\bu_{ht}^{*},\bphi_h)+b(\bu_H,\bu_{ht}^{*},\bphi_h)
 +b(\bu_{ht}^{*},\bu_H,\bphi_h)=(\f_t,\bphi_h)
\nonumber\\
&-\big(b(\bu_{Ht},\bu_{h}^{*},\bphi_h)+b(\bu_h^{*},\bu_{Ht},\bphi_h)\big)+\big(b(\bu_H,\bu_{Ht},\bphi_h)
+b(\bu_{Ht},\bu_{H},\bphi_h)\big).
\end{align}
Substitute $\bphi_h=e^{2\alpha t}\bu_{ht}^{*}$ in (\ref{e11*}) and use (\ref{tw11}) to obtain
\begin{align}\label{e12}
 \frac{e^{2\alpha t}}{2}\frac{d}{dt}\|\bu_{ht}^{*}\|^2+\gamma e^{2\alpha t}\|\nabla\bu_{ht}^{*}\|^2\leq e^{2\alpha t}(\f_t,\bu_{ht}^{*})+I,~ \rm{say}.
\end{align}
A use of Lemma \ref{nlt} yields
\begin{align}\label{e13}
|I|&\leq e^{2\alpha t}|b(\bu_{Ht},\bu_h^{*},\bu_{ht}^{*})|+|b(\bu_{h}^{*},\bu_{Ht},\bu_{ht}^{*})|
+|b(\bu_{Ht},\bu_H,\bu_{ht}^{*})|+|b(\bu_H,\bu_{Ht},\bu_{ht}^{*})|\nonumber\\
&\leq C e^{2\alpha t} \|\bu_{Ht}\|(\|\tilde\Delta_h\bu_{h}^{*}\|+\|\tilde\Delta_H\bu_{H}\|)
\|\nabla\bu_{ht}^{*}\|.
\end{align}
Apply (\ref{2.1*}), (\ref{e13}), Cauchy-Schwarz's inequality and Young's inequality in (\ref{e12}) to arrive at
\begin{align}
 \frac{d}{dt}e^{2\alpha t}\|\bu_{ht}^{*}\|^2+\left(\gamma-\frac{2\alpha}{\lambda_1}\right) 
 e^{2\alpha t}\|\nabla\bu_{ht}^{*}\|^2\leq 
C e^{2\alpha t}\left(\|\f_t\|^2
+\|\bu_{Ht}\|^2(\|\tilde\Delta_h\bu_h^{*}\|^2+\|\tilde\Delta_H\bu_H\|^2)\right).\nonumber
\end{align}
An integration with respect to time, a use of assumption {\bf (A2)}, (\ref{e19}) and Lemmas \ref{1l1*}, \ref{1l1} leads to
\begin{align}\label{e15}
 \|\bu_{ht}^{*}(t)\|^2+e^{-2\alpha t}\displaystyle{\int_0^t} 
 e^{2\alpha s}\|\nabla\bu_{ht}^{*}(s)\|^2 ds\leq C.
\end{align}
\noindent
Next, choose $\bphi_h=-e^{2\alpha t}\tilde\Delta_h^{-1}\bu_{htt}^{*}$ in (\ref{e11*}) and use Lemma \ref{nlt} to arrive at
\begin{align}
\label{e20*}
|b(\bu_{Ht},\bu_H,\tilde\Delta_h^{-1}\bu_{htt}^{*})|+|b(\bu_{H},\bu_{Ht},\tilde\Delta_h^{-1}\bu_{htt}^{*})|&
\leq C\|\nabla\bu_H\|\|\nabla\bu_{Ht}\|\|\bu_{htt}^{*}\|_{-1},\nonumber\\
|b(\bu_{Ht},\bu_h^{*},\tilde\Delta_h^{-1}\bu_{htt}^{*})+|b(\bu_h^{*},\bu_{Ht},\tilde\Delta_h^{-1}\bu_{htt}^{*})|
&\leq C\|\nabla\bu_{Ht}\|\|\nabla\bu_h^{*}\|\|\bu_{htt}^{*}\|_{-1},\\
 |b(\bu_{H},\bu_{ht}^{*},\tilde\Delta_h^{-1}\bu_{htt}^{*})+|b(\bu_{ht}^{*},\bu_{H},\tilde\Delta_h^{-1}\bu_{htt}^{*})|
&\leq C\|\nabla\bu_{ht}^{*}\|\|\nabla\bu_{H}\|\|\bu_{htt}^{*}\|_{-1}.\nonumber
\end{align}
Integrate with respect to time from $0$ to $t$, use (\ref{e19}), (\ref{e20*}) and Lemma \ref{1l1} to obtain
\begin{align}
 e^{2\alpha t}\|\bu_{ht}^{*}(t)\|^2+ \int_0^t e^{2\alpha s}&\|\bu_{htt}^{*}(s)\|_{-1}^2ds \leq 
(\|\tilde\Delta_H\bu_{0H}\|^2+\|\tilde\Delta_h\bu_{0h}^{*}\|^2)+C\bigg(\int_0^te^{2\alpha s}\|\bu_{ht}^{*}(s)\|^2ds\nonumber\\
&+\int_0^t e^{2\alpha s}(\|\nabla\bu_{Ht}(s)\|^2+\|\nabla\bu_{ht}^{*}(s)\|^2)ds 
+\int_0^te^{2\alpha s}\|\f_t(s)\|^2ds\bigg).\nonumber
\end{align}
A use of (\ref{e16}), (\ref{e15}), {\bf (A2)} and Lemma \ref{1l1*} concludes the proof of Lemma \ref{1l2}.\hfill{$\Box$}
\begin{ldf}\label{2bdlm}{\it (a priori bounds for $\bU_H^n$)}
With $\alpha>0$, choose $k_0$ small so that for $0<k\leq k_0$,
\begin{align}\label{sbd11}
1+\left(\frac{\gamma\lambda_1}{2}\right)k\geq e^{\alpha k}.
 \end{align}
Further, let $\bU_H^0=P_H\bu_{0H}$. Then, discrete solution $\bU_H^n$, $n\geq 1$ of 
(\ref{1bd1}) satisfies the following estimates: 
 \begin{align}
&\|\bU_H^n\|^2+ e^{-2\alpha t_{n}}~k \displaystyle{\sum_{i=1}^{n}}
 e^{2\alpha t_i}\|\nabla \bU_H^i\|^2
\leq C(\gamma,\nu,\alpha,\lambda_1) \left(e^{-2\alpha t_n}\|\bU_H^0\|^2+\|\f\|^2_{\infty}\right),\nonumber\\
&\|\nabla\bU_H^n\|^2+ e^{-2\alpha t_{n}}~k \displaystyle{\sum_{i=1}^{n}}
 e^{2\alpha t_i}\|\tilde\Delta_H\bU_H^i\|^2
\leq C(\gamma,\nu,\alpha,\lambda_1) \left(e^{-2\alpha t_n}\|\nabla\bU_H^0\|^2+\|\f\|^2_{\infty}\right),\nonumber
 \end{align}
where $\|\f\|_{\infty}=\|\f\|_{L^{\infty}(\bL^2)}$.
 \hfill{$\Box$}
 \end{ldf}
\begin{ldf}\label{tsbd1}{\it (estimates for $\e_H^n$)}
Let the assumptions of Lemma \ref{2bdlm} be satisfied. Also, let $u_H(t)$ be a solution of (\ref{1a1}) and 
$\e_H^n=\bU_H^n-\bu_H^n$, for $n\geq 1$. Then, for 
some positive constant $K_T$, that depends on $T$, there holds
\begin{align}
 \|{\bf e}_H^n\|^2+  k e^{-2\alpha t_n} \sum_{i=1}^{n} e^{2\alpha t_i} 
\|\nabla \e_H^i\|^2 \leq  K_T  k^2.\nonumber
\end{align}
\end{ldf}
\begin{ldf}\label{bdlm1} {\it (a priori bounds for $\bU^n$)}
Under the hypotheses of Lemma \ref{2bdlm}, the discrete solution $\bU^n$, $n \geq 1$ of (\ref{2bd1}) satisfies
 \begin{align}
 &\|\bU^n\|^2+  e^{-2\alpha t_{n}}~k \displaystyle{\sum_{i=1}^{n}}
 e^{2\alpha t_i}\|\nabla\bU^i\|^2
\leq C(\gamma,\nu,\alpha,\lambda_1,T) (e^{-2\alpha t_n}\|\bU^0\|^2+\|\f\|^2_{\infty}).\nonumber\\
&\|\nabla\bU^n\|^2+  e^{-2\alpha t_{n}}~k \displaystyle{\sum_{i=1}^{n}}
 e^{2\alpha t_i}\|\tilde\Delta_h\bU^i\|^2
\leq C(\gamma,\nu,\alpha,\lambda_1,T) (e^{-2\alpha t_n}\|\nabla\bU^0\|^2+\|\f\|^2_{\infty}).\nonumber
 \end{align}
 \end{ldf} 
\noindent
{\it Proof.} For $n=i$, multiply (\ref{2bd1}) by $e^{\alpha t_i}$, use $e^{\alpha t_i}\bar\partial_t\bU^i=
e^{\alpha k}\bar\partial_t\hat\bU^i-\left(\frac{e^{\alpha k}-1}{k}\right)\hat\bU^i$ and 
divide the resulting equation by $e^{\alpha k}$ to obtain
\begin{align}\label{be1}
&(\bar\partial_t\hat\bU^i,\bphi_h)-\left(\frac{1-e^{-\alpha k}}{k}\right)(\hat\bU^i,\bphi_h)
+\nu\, e^{\alpha t_i}e^{-\alpha k}(a(\bU^i,\bphi_h)+b(\bu_H^i,\bU^i,\bphi_h)+b(\bU^i,\bu_H^i,\bphi_h))\nonumber\\
&=e^{-\alpha k}({\hat\f^i},\bphi_h)-e^{-\alpha {t_{i+1}}}b(\hat\e_H^i,\hat\bU^i,\bphi_h)
-e^{-\alpha {t_{i+1}}}b(\hat\bU^i,\hat\e_H^i,\bphi_h)
+e^{-\alpha {t_{i+1}}}
b(\hat\bU_H^i,\hat\bU_H^i,\bphi_h).
\end{align}
 Observe that
\begin{align}\label{be2}
 (\bar\partial_t\bphi^i,\bphi^i)=\frac{1}{2k}\left(\bphi^i-\bphi^{i-1}\right)\geq 
\frac{1}{2}\bar\partial_t\|\bphi^i\|^2.
\end{align}
Substitute $\bphi_h=\hat\bU^i$ in (\ref{be1}), use (\ref{3.1*}) and (\ref{be2}) to arrive at
\begin{align}\label{be3}
 \frac{1}{2}&\bar\partial_t\|\hat\bU^i\|^2+\left(\gamma e^{-\alpha k}-
 \left(\frac{1-e^{-\alpha k}}{k}\right)\lambda_1^{-1}\right)
\|\nabla\hat\bU^i\|^2\nonumber\\
&\leq e^{-\alpha k}({\hat\f^i},\hat\bU^i)
-e^{-\alpha t_{i+1}}b(\hat\bU^i,\hat\e_H^i,\hat\bU^i)+e^{-\alpha t_{i+1}}
b(\hat\bU_H^i,\hat\bU_H^i,\hat\bU^i).
\end{align}
Applying (\ref{2.1*}) and Cauchy-Schwarz's inequality, the first term on the right hand side of (\ref{be3}) can be bounded as 
\begin{align}\label{be4}
 |e^{-\alpha k}({\hat\f^i},\hat\bU^i)|\leq Ce^{-\alpha k}\|\hat\f^i\|\|\hat\bU^i\|\leq C(\lambda_1,\epsilon)e^{-\alpha k}
\|\hat\f^i\|^2+\epsilon e^{-\alpha k}\|\nabla\hat\bU^i\|^2.
\end{align}
A use of Lemma \ref{nlt} with Young's inequality yields
\begin{align}\label{be5}
 e^{-\alpha t_{i+1}}(|b(\hat\bU_H^i,\hat\bU_H^i,\hat\bU^i)|+|b(\hat\bU^i,\hat\e_H^i,\hat\bU^i)|)
& \leq C(\epsilon)e^{-2\alpha t_i}e^{-\alpha k}(\|\nabla\hat\bU_H^i\|^4+\|\nabla\hat\e_H^i\|^2\|\hat\bU^i\|^2)\nonumber\\
&+\epsilon e^{-\alpha k}\|\nabla\hat\bU^i\|^2.
\end{align}
Apply (\ref{be4})-(\ref{be5}) with $\epsilon=\gamma/2$ in (\ref{be3}) to obtain
\begin{align}\label{be7}
  &\bar\partial_t\|\hat\bU^i\|^2+\left(\gamma\,e^{-\alpha k}-2\left(\frac{1-e^{-\alpha k}}{k}\right)\lambda_1^{-1}\right)
\|\nabla\hat\bU^i\|^2\nonumber\\
&\leq Ce^{-\alpha k}\left(\|{\hat\f^i}\|^2+e^{-2\alpha t_i}\|\nabla\hat\bU_H^i\|^4
+e^{-2\alpha t_i}\|\nabla\hat\e_H^i\|^2\|\hat\bU^i\|^2\right).
\end{align}
We choose $k_0>0$, such that $ 1+\left(\frac{\gamma\lambda_1}{2}\right)k\geq e^{\alpha k}$. 
This guarantees that $\gamma e^{-\alpha k}-2\left(\frac{1-e^{-\alpha k}}{k}\right)\lambda_1^{-1}\geq 0$. 
Multiply (\ref{be7}) by $k$ and then sum over $i=1$ to $n$  to obtain
\begin{align}
\|\hat\bU^n\|^2+\bigg(\gamma e^{-\alpha k}-2\bigg(\frac{1-e^{-\alpha k}}{k}\bigg)&\lambda_1^{-1}\bigg)k
\displaystyle{\sum_{i=1}^{n}}\|\nabla\hat\bU^i\|^2\leq \|\bU^0\|^2+C\bigg(ke^{-\alpha k}\|{\f}\|_{\infty}^2
\displaystyle{\sum_{i=1}^{n}}e^{2\alpha t_i}\nonumber\\
&+k\displaystyle{\sum_{i=1}^{n}}e^{-2\alpha t_i}\|\nabla\hat\bU_H^i\|^4
+k\displaystyle{\sum_{i=1}^{n}}e^{-2\alpha t_i}\|\nabla\hat\e_H^i\|^2\|\hat\bU^i\|^2\bigg).\nonumber
\end{align}
An application of Gronwall's lemma with Lemmas \ref{2bdlm} and \ref{tsbd1} leads to the desired result.\\
For $n=i$, multiply (\ref{2bd1}) by $e^{2\alpha t_i}$ and choose $\bphi_h=-\tilde \Delta_h\hat\bU^i$ to arrive at
\begin{align}\label{7n}
 \frac{1}{2}&\bar\partial_t\|\nabla\hat\bU^i\|^2+\nu e^{-\alpha k}\|\tilde\Delta_h\hat\bU^i\|^2\leq 
 -e^{-\alpha k}(\hat\f^i,\tilde\Delta_h\hat\bU^i)
 -\left(\frac{1-e^{-\alpha k}}{k}\right)(\hat\bU^i,\tilde\Delta_h\hat\bU^i)\nonumber\\
 &+e^{-\alpha t_{i+1}}(b(\hat\bU_H^i,\hat\bU^i,\tilde\Delta_h\hat\bU^i)+b(\hat\bU^i,\hat\bU_H^i,\tilde\Delta_h\hat\bU^i))
 -e^{-\alpha t_{i+1}}b(\hat\bU_H^i,\hat\bU_H^i,\tilde\Delta_h\hat\bU^i).
\end{align}
A use of Lemma \ref{nlt} leads to 
\begin{align}\label{7n1}
 |b(\hat\bU_H^i,\hat\bU^i,\tilde\Delta_h\hat\bU^i)|&+|b(\hat\bU^i,\hat\bU_H^i,\tilde\Delta_h\hat\bU^i)|
 +|b(\hat\bU_H^i,\hat\bU_H^i,\tilde\Delta_h\hat\bU^i)|\nonumber\\
 &\leq C (\|\tilde\Delta_H\hat\bU_H^i\|\|\nabla\hat\bU^{i}\|+
 \|\nabla\hat\bU_H^i\|\|\tilde\Delta_H\hat\bU_H^i\|)\|\tilde\Delta_h\hat\bU^i\|.
\end{align}
Multiply (\ref{7n}) by $k$ and then sum over $i=1$ to $n$ and use (\ref{2.1*}), (\ref{7n1}) to arrive at
\begin{align}
 \|\nabla\hat\bU^n\|^2&+\nu e^{-\alpha k}k\displaystyle{\sum_{i=1}^{n}}\|\tilde\Delta_h\hat\bU^i\|^2\leq 
 \|\nabla\bU^0\|^2+C(\lambda_1)k\displaystyle{\sum_{i=1}^{n}}e^{-\alpha k}(\|\hat\f^i\|^2+\|\nabla\hat\bU^i\|^2 )\nonumber\\
 &+e^{-\alpha k}k\displaystyle{\sum_{i=1}^{n}}e^{-2\alpha t_{i}}(\|\tilde\Delta_H\hat\bU_H^i\|^2\|\nabla\hat\bU^{i}\|^2+
 \|\nabla\hat\bU_H^i\|^2\|\tilde\Delta_H\hat\bU_H^i\|^2).
\end{align}
An application of Gronwall's lemma with {\bf(A2)}, Lemmas \ref{2bdlm} and \ref{bdlm1} concludes the proof.\hfill{$\Box$}
\begin{ldf}\label{bdlm2} {\it (a priori bounds for $\bU_h^n$)}
Under the hypotheses of Lemma \ref{2bdlm}, the discrete solution $\bU_h^n$, $n \geq 1$ of (\ref{3bd1}) satisfies
 \begin{align}
 &\|\bU_h^n\|^2+  e^{-2\alpha t_{n}}~k \displaystyle{\sum_{i=1}^{n}}
 e^{2\alpha t_i}\|\nabla\bU_h^i\|^2
\leq C(\gamma,\nu,\alpha,\lambda_1,T)(e^{-2\alpha t_n}\|\bU^0\|^2+\|\f\|^2_{\infty}),\nonumber\\
&\|\nabla\bU_h^n\|^2+  e^{-2\alpha t_{n}}~k \displaystyle{\sum_{i=1}^{n}}
 e^{2\alpha t_i}\|\tilde\Delta_h\bU_h^i\|^2
\leq C(\gamma,\nu,\alpha,\lambda_1,T) (e^{-2\alpha t_n}\|\nabla\bU^0\|^2+\|\f\|^2_{\infty}).\nonumber
 \end{align}
 \end{ldf} 
 \noindent
 {\it Proof.} For $n=i$, multiply (\ref{3bd1}) by $e^{\alpha t_i}$, substitute $\bphi_h=\hat\bU_h^i$, 
 use (\ref{3.1*}) and (\ref{be2}) to arrive at
\begin{align}\label{ee19}
 \frac{1}{2}&\bar\partial_t\|\hat\bU_h^i\|^2+\left(\gamma e^{-\alpha k}
 -\left(\frac{1-e^{-\alpha k}}{k}\right)\lambda_1^{-1}\right)
\|\nabla\hat\bU_h^i\|^2\leq e^{-\alpha k}({\hat\f^i},\hat\bU_h^i)\nonumber\\
&-e^{-\alpha t_{i+1}}b(\hat\bU_h^i,\hat\e_H^i,\hat\bU_h^i)+e^{-\alpha t_{i+1}}
\left(b(\hat\bU_H^i,\hat\bU^{i},\hat\bU_h^i)+b(\hat\bU^{i},\hat\bU_H^i-\hat\bU^{i},\hat\bU^i_h)\right).
\end{align}
The first two terms on the right hand side of (\ref{ee19}) can be tackled similar to (\ref{be4})-(\ref{be5}). 
To bound the third term, use Lemmas \ref{nlt}, \ref{2bdlm}, \ref{bdlm1} and Young's inequality and arrive at 
\begin{align}\label{ee20}
e^{-\alpha t_{i+1}} |b(\hat\bU_H^i,\hat\bU^{i},\hat\bU_h^i)+b(\hat\bU^{i},\hat\bU_H^i-\hat\bU^{i},\bU_h^i)|\leq& 
C e^{-2\alpha t_{i}}e^{-\alpha k}(\|\nabla \hat\bU_H^i\|^2+\|\nabla \hat\bU^i\|^2)\nonumber\\
&+\epsilon e^{-\alpha k}
\|\nabla\hat\bU_h^i\|^2.
\end{align}
A use of (\ref{be4})-(\ref{be5}) and (\ref{ee20}) in (\ref{ee19}) yields
 \begin{align}\label{ee21}
  &\bar\partial_t\|\hat\bU_h^i\|^2+\left(\gamma e^{-\alpha k}-2\left(\frac{1-e^{-\alpha k}}{k}\right)
\lambda_1^{-1}\right)
\|\nabla\hat\bU_h^i\|^2\leq Ce^{-\alpha k}\times\nonumber\\
&\left(\|{\hat\f^i}\|^2+e^{-2\alpha t_i}(\|\nabla\hat\bU_H^i\|^2+\|\nabla\hat\bU^i\|^2)
+e^{-2\alpha t_i}\|\nabla\hat\e_H^i\|^2\|\hat\bU_h^i\|^2\right).
\end{align}
Multiply (\ref{ee21}) by $k$, sum over $i=1$ to $n$ and use Lemmas \ref{2bdlm}-\ref{bdlm1} to complete the proof.
\hfill{$\Box$}
\subsection{{\it A Priori} Error Estimates}
Consider (\ref{1a1})-(\ref{3a1}) at $t=t_n$ and subtract the resulting equations from (\ref{1bd1})-(\ref{3bd1}), 
respectively, to arrive at the following error equations: \\
{\it {\bf Step 1.}} for all $\bphi_H\in \bJ_H$
\begin{eqnarray}
&&\hspace{-1cm}{\footnotesize(\bar \partial_t\e_H^n,\bphi_H)+\nu\,a(\e_H^n,\bphi_H)+b(\bu_H^n,\e_H^n,\bphi_H)
+b(\e_H^n,\bu_H^n,\bphi_H)=(\bs_H^n,\bphi_H)
+\Lambda_H(\bphi_H)},\label{sbd3n}
\end{eqnarray}
where {\footnotesize$\Lambda_H(\bphi_H)=b(\bu_H^n,\e_H^n,\bphi_H)-b(\bU_H^n,\e_H^n,\bphi_H)$ and 
$\bs_{H}^n=\bu_{Ht}^{n}-\bar\partial_t\bu_H^n$.} \\
\noindent
 {\it {\bf Step 2.}} for all $\bphi_h\in \bJ_h$
 \begin{eqnarray}
 \hspace{-1cm}{\footnotesize(\bar\partial_t\e^n,\bphi_h)+\nu\,a(\e^n,\bphi_h)+b(\bu_H^n,\e^n,\bphi_h)
+b(\e^n,\bu_H^n,\bphi_h)=
(\bs^{n},\bphi_h)+\Lambda^{*}(\bphi_h)},\label{sbd4n}
 \end{eqnarray}
where {\footnotesize $\bs^{n}= \bu_{ht}^{*n}-\bar\partial_t \bu_h^{*n}$,}
 { \footnotesize$\Lambda^{*}(\bphi_h)=\Lambda_1(\bphi_h)+\Lambda_2(\bphi_h)+\Lambda_3(\bphi_h)$}
 with \begin{eqnarray}\label{sbd5}
 \left.
 \begin{array}{rcl}
 &&{\footnotesize \Lambda_1(\bphi_h)=-b(\e_H^n,\bU^n,\bphi_h),}\\ 
 &&{\footnotesize\Lambda_2(\bphi_h)=-b(\bU^n,\e_H^n,\bphi_h),} \\ 
  &&{\footnotesize\Lambda_3(\bphi_h)=b(\bU_H^n,\bU_H^n,\bphi_h)-b(\bu_H^n,\bu_H^n,\bphi_h)}.\\
 \end{array}
 \right\}
 \end{eqnarray}
 {\it {\bf Step 3.}} for all $\bphi_h\in \bJ_h$
 \begin{eqnarray}
 &&\hspace{-1cm} {\footnotesize(\bar\partial_t\e_h^n,\bphi_h)+\nu\,a(\e_h^n,\bphi_h)+b(\bu_H^n,\e_h^n,\bphi_h)
+b(\e_h^n,\bu_H^n,\bphi_h)
 =(\bs_{h}^n,\bphi_h)+\Lambda_h(\bphi_h),}\label{sbd7n}
 \end{eqnarray}
 where {\footnotesize $\bs_h^n= \bu_{ht}^{n}-\bar\partial_t \bu_h^{n}$,} 
 { \footnotesize$\Lambda_h(\bphi_h)= \Lambda^1_h(\bphi_h)+\Lambda^2_h(\bphi_h)+\Lambda^3_h(\bphi_h)+\Lambda^4_h(\bphi_h)$} 
with 
 \begin{eqnarray}\label{sbd8}
  \left.
 \begin{array}{rcl}
 &&{\footnotesize\Lambda^1_h(\bphi_h)=-b(\e_H^n,\bU_h^n,\bphi_h),}\\
 &&{\footnotesize\Lambda^2_h(\bphi_h)=-b(\bU_h^n,\e_H^n,\bphi_h),}\\
  &&{\footnotesize\Lambda^3_h(\bphi_h)=b(\bU_H^n,\bU^n,\bphi_h)-b(\bu_H^n,\bu_h^{*n},\bphi_h),}\\
 &&{\footnotesize\Lambda^4_h(\bphi_h)=b(\bU^n,\bU_H^n-\bU^n,\bphi_h)-b(\bu_h^{*n},\bu_H^n-\bu_h^{*n},\bphi_h).}
 \end{array}
 \right\}
 \end{eqnarray}
\noindent
The main result of this section is stated as:
\begin{tdf}\label{Tbd1} {\it (fully discrete error estimates)}
   Under the assumptions of Theorem \ref{T5} and Lemma \ref{2bdlm}, the following hold true:  
   \begin{eqnarray}
   && \|\bu(t_n)-\bU_h^n\|\leq C (h^2+H^{4-2\delta}+k),\qquad \|\nabla(\bu(t_n)-\bU_h^n)\|\leq C (h+H^{4-\delta}+k),\nonumber\\
  &&\|p(t_n)-P_h^n\|\leq C (h+H^{4-\delta}+k^{1/2}),\nonumber
 \end{eqnarray}
 where $\delta>0$ is arbitrarily small.
\end{tdf}
Below, we prove a lemma which will be used subsequently.
\begin{ldf}\label{tsbd2}
  Assume that \rm{\bf (A1)}-\rm{\bf (A2)} and \rm{\bf (B1)}-\rm{\bf (B2)} hold true. Let 
 for some fixed $h$, $\bu_h^{*}$ satisfies (\ref{2a1}). 
 Then, there is a positive constant $K_T$ that depends on $T$ such that
 \begin{align}
  \|\e^i\|^2+ke^{-2\alpha t_n}\displaystyle{\sum_{i=1}^n}e^{2\alpha t_i}\|\nabla\e^i\|^2\leq K_T k^2.\nonumber
 \end{align}
 \end{ldf}
\noindent
{\it Proof.} For $n=i$, substitute $\bphi_h=\e^i$ in (\ref{sbd4n}) and use (\ref{be2}) to obtain
\begin{align}\label{ee1}
 \bar\partial_t\|\e^i\|^2+2\gamma\|\nabla\e^i\|^2\leq 2(\bs^{i},\e^i)+2\Lambda^{*}(\e^i).
\end{align}
Multiply (\ref{ee1}) by $e^{2\alpha ik}$ and sum over $i=1$ to $n$, where $T=nk$. Use the fact  
\begin{align}\label{ee2}
 \displaystyle{\sum_{i=1}^{n}} k e^{2\alpha ik}\bar\partial_t\|\e^i\|^2&= \displaystyle{\sum_{i=1}^{n}}
e^{2\alpha ik}(\|\e^i\|^2-\|\e^{i-1}\|^2)\nonumber\\
&=e^{2\alpha nk}\|\e^n\|^2-\displaystyle{\sum_{i=1}^{n-1}}e^{2\alpha {i}k}(e^{2\alpha k}-1)\|\e^i\|^2 
\end{align}
to arrive at
\begin{align}\label{ee3}
 e^{2\alpha nk}\|\e^n\|^2+2k\gamma&\displaystyle{\sum_{i=1}^{n}}e^{2\alpha ik}\|\nabla\e^i\|^2\leq 
\displaystyle{\sum_{i=1}^{n-1}}e^{2\alpha {i}k}(e^{2\alpha k}-1)\|\e^i\|^2\nonumber\\
&+2k\displaystyle{\sum_{i=1}^{n}}e^{2\alpha {i}k}
(\bs^{i},\e^i)+2k \displaystyle{\sum_{i=1}^{n}}e^{2\alpha {i}k}\Lambda^{*}(\e^i).
\end{align}
A use of Taylor's series expansion in the interval $(t_{i-1},t_i)$ with use of Cauchy-Schwarz's and Young's inequalities yields
\begin{align}\label{ee4}
|2(\bs^{i},\e^i)|&\leq \frac{2}{k}\displaystyle{\int_{t_{i-1}}^{t_i}}(t-t_{i-1})\|\bu_{htt}^{*}\|_{-1}dt\|\nabla\e^i\|\nonumber\\
 & \leq K\,k^{1/2}\left\{ \displaystyle{\int_{t_{i-1}}^{t_i}}\|\bu_{htt}^{*}\|^2_{-1}dt\right\}^{1/2}\|\nabla\e^i\|. 
\end{align}
From Lemma \ref{1l2}, observe that 
\begin{align}\label{ee5}
 \displaystyle{\sum_{i=1}^{n}}e^{2\alpha {i}k}\displaystyle{\int_{t_{i-1}}^{t_i}}\|\bu_{htt}^{*}\|^2_{-1}dt
&=\displaystyle{\sum_{i=1}^{n}}\displaystyle{\int_{t_{i-1}}^{t_i}}e^{2\alpha (t_i-t)}e^{2\alpha t}\|\bu_{htt}^{*}\|^2_{-1}dt\nonumber\\
&\leq e^{2\alpha k}\displaystyle{\int_{0}^{t_n}}e^{2\alpha t}\|\bu_{htt}^{*}\|^2_{-1}dt\leq Ke^{2\alpha (n+1)k}.
\end{align}
Apply (\ref{sbd5}) to obtain
\begin{align}\label{ee6}
 |2k \displaystyle{\sum_{i=1}^{n}}e^{2\alpha {i}k}\Lambda^{*}(\e^i)|\leq 2k \displaystyle{\sum_{i=1}^{n}}e^{2\alpha {i}k}
(|\Lambda_1(\e^i)|+|\Lambda_2(\e^i)|+|\Lambda_3(\e^i)|).
\end{align}
An application of Lemma \ref{nlt} with Young's inequality and Lemmas \ref{tsbd1} and \ref{bdlm1} leads to
\begin{align}\label{ee7}
 2k\displaystyle{\sum_{i=1}^{n}}e^{2\alpha {i}k}(|\Lambda_1(\e^i)|+|\Lambda_2(\e^i)|)
&\leq C k\displaystyle{\sum_{i=1}^{n}}e^{2\alpha {i}k}\|\nabla \e_H^i\|\|\nabla\bU^i\|\|\nabla\e^i\|\nonumber\\
&\leq C(\epsilon) k\displaystyle{\sum_{i=1}^{n}}e^{2\alpha {i}k}\|\nabla \e_H^i\|^2\|\nabla\bU^i\|^2+\epsilon k
\displaystyle{\sum_{i=1}^{n}}e^{2\alpha {i}k}\|\nabla\e^i\|^2\nonumber\\
& \leq C(T,\epsilon)k^2 e^{2\alpha nk}+\epsilon k\displaystyle{\sum_{i=1}^{n}}e^{2\alpha {i}k}\|\nabla\e^i\|^2.
\end{align}
A use of boundedness of $\|\nabla\bu_H\|\leq C$ and Lemmas \ref{nlt} and \ref{2bdlm} to obtain 
\begin{align}
 |\Lambda_3(\e^i)|&= |b(\e_H^i, \bu_H^i,\e^i)
+b(\bU_H^i, \e_H^i,\e^i)|\nonumber\\
&\leq C (\|\nabla\bU_H^i\|+\|\nabla\bu_H^i\|)\|\nabla\e_H^i\|\|\nabla\e^i\|
\leq C \|\nabla\e_H^i\|\|\nabla\e^i\|.\nonumber
\end{align}
Now, a use of Young's inequality and Lemma \ref{tsbd1} yields
\begin{align}\label{ee9}
 |2k \displaystyle{\sum_{i=1}^{n}}e^{2\alpha {i}k}\Lambda_3(\e^i)|&\leq C(\epsilon)
 k\displaystyle{\sum_{i=1}^{n}}e^{2\alpha {i}k}
\|\nabla\e_H^i\|^2+\epsilon k \displaystyle{\sum_{i=1}^{n}}e^{2\alpha {i}k} \|\nabla\e^i\|^2\nonumber\\
&\leq C(T,\epsilon)k^2 e^{2\alpha nk}+\epsilon k \displaystyle{\sum_{i=1}^{n}}e^{2\alpha {i}k} \|\nabla\e^i\|^2.
\end{align}
With the help of (\ref{ee4})-(\ref{ee9}), (\ref{ee3}) can be written as
\begin{align}
 e^{2\alpha nk}\|\e^n\|^2+k\displaystyle{\sum_{i=1}^{n}}e^{2\alpha {i}k}\|\nabla\e^i\|^2\leq 
 Ck^2(e^{2\alpha(n+1)k}+e^{2\alpha nk})+C k \displaystyle{\sum_{i=1}^{n-1}}e^{2\alpha {i}k}\|\e^i\|^2.\nonumber
\end{align}
 A use of discrete Gronwall's Lemma leads to
\begin{align}
 \|\e^n\|^2+ke^{-2\alpha nk}\displaystyle{\sum_{i=1}^{n}}e^{2\alpha ik}\|\nabla\e^i\|^2\leq K_T k^2\nonumber 
\end{align}
and this completes the rest of the proof.  \hfill{$\Box$}

Now, substitute $\bphi_h=(-\tilde \Delta_h)^{-1}\bar\partial_t \e^n$ in (\ref{sbd4n}) and use Lemma \ref{nlt} to arrive at
\begin{align}
\|\bar\partial_t\e^n\|^2_{-1}&\leq C(\nu)(\|\nabla\e^n\|+\|\bs^{n}\|+\|\nabla\bu_H^n\|\|\nabla\e^n\|\nonumber\\
&+(\|\nabla\bU_H^n\|+\|\nabla\bu_H^n\|+\|\nabla\bU^n\|)\|\nabla\e_H^n\|)\|\bar\partial_t\e^n\|_{-1}.
\end{align}
An application of (\ref{ee4}) and Lemmas \ref{1l1*}, \ref{1l2}, \ref{2bdlm}, \ref{tsbd1}, \ref{bdlm1}, \ref{tsbd2} leads to
\begin{align}\label{ee4*}
\|\bar\partial_t\e^n\|^2_{-1}&\leq Ck.
\end{align}
Next, to derive pressure error estimates, we consider the equivalent form of semidiscrete approximations (\ref{2a1}) as:
find $(\bu_h^{*}(t),p_h^{*}(t)) \in {\bf H}_h \times L_h$ such that $ \bu_h^{*}(0)= \bu_{0h} $ and for $t>0$,
 \begin{eqnarray}\label{2a1*}
 \left.
 \begin{array}{rcl}
&&(\bu_{ht}^{*},\bphi_h)+\nu a(\bu_h^{*}, \bphi_h) 
+b(\bu_h^{*},\bu_H,\bphi_h)\\
&&\hspace{.2cm}+\,b(\bu_H,\bu_h^{*},\bphi_h)=b(\bu_H,\bu_H,\bphi_h)
+(p_h^{*}, \nabla \cdot {\bphi}_h)\;\;\; \forall {\bphi}_h \in {\bf H}_h,\\
 \hspace{-.1cm}&&(\nabla \cdot \bu_h^{*}, \chi_h) = 0 \;\;\; \forall \chi_h \in L_h.
\end{array}
\right\}
 \end{eqnarray}
The equivalent form of fully discrete approximations (\ref{2bd1}) is as follows: $\forall (\bphi_h,\chi_h) \in \bH_h\times 
 L_h$, seek a sequence of functions  ${(\bU^n, P^n)}_{n\geq1}\in \bH_h\times L_h$ as solutions of the following equations:
\begin{eqnarray} \label{a58}
 \left.
 \begin{array}{rcl}
&&(\bar\partial_t \bU^n,\bphi_h)+\nu\, a(\bU^n, \bphi_h) 
+b(\bU^n,\bU_H^n,\bphi_h)\\
&&+b(\bU^n_H,\bU^n,\bphi_h)=b(\bU_H^n,\bU_H^n,\bphi_h)+(P^n,\nabla\cdot\bphi_h),\\
&&(\nabla\cdot\bU^n,\chi_h)=0.
 \end{array}
 \right\}
 \end{eqnarray}
Subtract (\ref{2a1*}) from (\ref{a58}) and write $\brho^n=P^n- p_h^{*n}$ to obtain 
\begin{align}
 (\brho^n,\nabla\cdot\bphi_h)&= (\bar\partial_t \e^n,\bphi_h)+\nu a(\e^n,\bphi_h)-\Lambda^{*}(\bphi_h)-(\bs^n,\bphi_h).\nonumber
\end{align}
A use of the Cauchy-Schwarz's inequality along with (\ref{ee4}), (\ref{ee4*}) and Lemmas \ref{1l1*}, \ref{1l2}, 
\ref{2bdlm}, \ref{tsbd1}, \ref{bdlm1}, \ref{tsbd2} yields
\begin{align}\label{a60}
 \|\brho^n\|\leq C(\kappa,\nu,\lambda_1,M) k^{1/2}.
\end{align}
A combination of (\ref{a60}) and Theorem \ref{tp1n} leads to the following pressure estimate.
\begin{align}
 \|p(t_n)-P^n\|\leq C (h+H^{3-\delta}+k^{1/2}).\nonumber
\end{align}

{\it Proof of Theorem \ref{Tbd1}.} Write $\bu(t_n)-\bU_h^n=(\bu(t_n)-\bu_h(t_n))-\e_h^n$. 
The estimate of $\bu(t_n)-\bu_h(t_n)$ is obtained in 
 Theorem \ref{T5}. Next, we proceed to derive the estimates for $\e_h^n$. For $n=i$, 
 substitute $\bphi_h=\e_h^i$ in (\ref{sbd7n}) and use (\ref{be2}). Multiply the resulting equation by 
$e^{2\alpha ik}$ and sum over $i=1$ to $n$ to obtain  
\begin{align}\label{ee11}
 e^{2\alpha nk}\|\e_h^n\|^2+2k\gamma &\displaystyle{\sum_{i=1}^{n}}e^{2\alpha ik}\|\nabla\e_h^i\|^2\leq 
\displaystyle{\sum_{i=1}^{n-1}}e^{2\alpha {i}k}(e^{2\alpha k}-1)\|\e_h^i\|^2\nonumber\\
&+2k\displaystyle{\sum_{i=1}^{n}}e^{2\alpha {i}k}
(\bs_h^{i},\e_h^i)+2 k \displaystyle{\sum_{i=1}^{n}}e^{2\alpha {i}k}\Lambda_h(\e_h^i).
\end{align}
The second term in the right hand side of (\ref{ee11}) can be bounded similar to (\ref{ee4})-(\ref{ee5}). Also, from (\ref{sbd8}) observe that
\begin{align}\label{ee12}
 2k \displaystyle{\sum_{i=1}^{n}}e^{2\alpha {i}k}|\Lambda_h(\e_h^i)|
 \leq 2k \displaystyle{\sum_{i=1}^{n}}e^{2\alpha {i}k}
\left(|\Lambda_h^1(\e_h^i)|+|\Lambda_h^2(\e_h^i)|+|\Lambda_h^3(\e_h^i)|+|\Lambda_h^4(\e_h^i)|\right).
\end{align}
 An application of Lemmas \ref{nlt}, \ref{tsbd1} and \ref{bdlm2} yields 
\begin{align}\label{ee13} 
 2k\displaystyle{\sum_{i=1}^{n}}e^{2\alpha {i}k}(|\Lambda_h^1(\e_h^i)|+|\Lambda_h^2(\e_h^i)|)
& \leq C(T,\epsilon)k^2e^{2\alpha nk}+k\epsilon\displaystyle{\sum_{i=1}^{n}}e^{2\alpha {i}k}\|\nabla\e_h^i\|^2.
\end{align}
Observe that
 \begin{align}\label{ee14}
  |\Lambda^3_h(\e_h^i)|\leq|b(\bU_H^i,\e^i,\e_h^i)|+|b(\e_H^i,\bu_h^{*i},\e_h^i)|.
 \end{align}
With the help of Lemmas \ref{nlt}, \ref{1l1}, \ref{2bdlm} and Young's inequality, it follows that
\begin{eqnarray}\label{ee15}
 2k \displaystyle{\sum_{i=1}^{n}}e^{2\alpha {i}k} |\Lambda^3_h(\e_h^i)|
\leq C(\epsilon)\displaystyle{\sum_{i=1}^{n}}e^{2\alpha {i}k}\big(\|\nabla\e^i\|^2+\|\nabla\e_H^i\|^2\big)
+\epsilon k \displaystyle{\sum_{i=1}^{n}}e^{2\alpha {i}k}\|\nabla\e_h^i\|^2.  
\end{eqnarray}
For the estimation of the fourth term on the right hand side of (\ref{ee12}), rewrite it as 
\begin{align}\label{ee16}
 |\Lambda^4_h(\e_h^i)|
=|b(\bU^i,\e_H^i,\e_h^i)-b(\bU^i,\e^i,\e_h^i)+b(\e^i,\bu_H^i-\bu_h^{*i},\e_h^i)|.					
\end{align}
Apply Lemma \ref{nlt}, \ref{1l1*}, \ref{1l1}, \ref{bdlm1} and Young's inequality to obtain
\begin{eqnarray}\label{ee17}
2k \displaystyle{\sum_{i=1}^{n}}e^{2\alpha {i}k} |\Lambda^4_h(\e_h^i)|
\leq  C(T,\epsilon) k \displaystyle{\sum_{i=1}^{n}}e^{2\alpha {i}k}
\big(\|\nabla \e^i\|^2+\|\nabla\e_H^i\|^2\big)+\epsilon k \displaystyle{\sum_{i=1}^{n}}e^{2\alpha {i}k}\|\nabla\e_h^i\|^2. 
\end{eqnarray}
A use of (\ref{ee13})-(\ref{ee17}) in (\ref{ee11}) leads to
\begin{align}\label{ee18}
 e^{2\alpha nk}\|\e_h^n\|^2+k\displaystyle{\sum_{i=1}^{n}}e^{2\alpha {i}k}\|\nabla\e_h^i\|^2\leq 
 Ck^2(e^{2\alpha(n+1)k}+e^{2\alpha nk})+C k \displaystyle{\sum_{i=1}^{n}}e^{2\alpha {i}k}\|\e_h^i\|^2.
\end{align}
 Use discrete Gronwall's Lemma to arrive at
\begin{align}\label{ee18*}
 \|\e_h^n\|^2+ke^{-2\alpha nk}\displaystyle{\sum_{i=1}^{n}}e^{2\alpha ik}\|\nabla\e_h^i\|^2\leq K_T k^2. 
\end{align}
A use of (\ref{ee18*}) along with Theorem \ref{T5} completes the proof of error estimates for velocity in Theorem \ref{Tbd1}.

Using the similar techniques as to arrive at (\ref{a60}) and Theorem \ref{T5}, the desired pressure estimate in 
Theorem \ref{Tbd1} can be obtained and this will conclude the proof of Theorem \ref{Tbd1}.\hfill{$\Box$}

\section{ Numerical Experiments}
 \setcounter{equation}{0}
 In this section, numerical results are presented to support theoretical results in Theorem \ref{Tbd1}. For 
space discretization, 
 $P_2$-$P_0$ mixed finite element space is used. We choose the domain $\Omega=(0,1)\times(0,1)$, time $t=[0,1]$, 
coefficients $\nu=1$ and $h=\mathcal{O}(H^2)$. Here, $N$ denotes the number of unknowns in the system. 
\begin{example}\label{ex1}
The right hand side function $f$ is chosen in such a way that the exact solution $(\bu,p)=((u_1,u_2),p)$ is $~~u_1  =   2e^{t} x^2(x-1)^2 y (y-1) (2y-1) $, $~~u_2  =  -2e^{t}y^2(y-1)^2 x (x-1) 
(2x-1)$, $~~~~p = y e^{t}$.
\end{example}
\noindent
 Table 1 gives the numerical errors and convergence rates obtained on successively 
refined meshes for backward Euler scheme with $k=\mathcal{O}(h^2)$ applied to two grid system 
(\ref{1a1})-(\ref{3a1}). The theoretical analysis provides a convergence rate of $\mathcal{O}(h^2)$ in $\bL^2$-norm, 
of $\mathcal{O}(h)$ in $\bH^1$-norm 
for velocity and of $\mathcal{O}(h)$ in $L^2$-norm for pressure with a choice of $k=\mathcal{O}(h)$. 
These results support the optimal theoretical convergence rates obtained in Theorem \ref{Tbd1} \\
\small{
\begin{table}[ht!]
 \centering
 \begin{tabular}{|l|l|l|l|l|l|l|l|l|l|}
 \hline
  \small{$N$}&  h &\small{$\|\bu(t_n)-\bU^n\|$}  &\small{Rate}& \small{$\|\bu(t_n)-\bU^n\|_{\bH^1}$} & \small{Rate} & \small{$\|p(t_n)-P^n\|$}&\small{ Rate}  \\
             &    &                                       &                 &            &        &                  &   \\
 \hline
 \hline      & 1/4 &0.009085  &           & 0.139927 &          &0.548331&   \\
 \hline 577  & 1/8 &0.002651  & 1.777183  & 0.075081 & 0.898156 &0.281244& 0.963220  \\
 \hline 2433 & 1/16&0.000713  & 1.893768  & 0.038833 & 0.951145 &0.142265& 0.983237  \\
 \hline 9986 & 1/32&0.000184  & 1.950443  & 0.019731 & 0.976861 &0.071518& 0.992191 \\
 \hline 40449& 1/64&0.000046  & 1.976824  & 0.009940 & 0.989066 &0.035856& 0.996088 \\
 \hline
 \end{tabular}
 \vspace{.1cm}
 \caption{ Errors and convergence rates for backward Euler method with $k=\mathcal{O}(h^2)$. }
 \end{table}
}
\begin{example}\label{ex2}
In this example, we choose the right hand side function
$f$ in such a way that the exact solution $(\bu,p)=((u_1,u_2),p)$ is:
 \begin{eqnarray} 
&&u_1  = t e^{-t^2} sin^2(3\pi x)~ sin(6\pi y),~~~~u_2  = -t e^{-t^2}sin^2(3 \pi y)~ sin(6 \pi x),\nonumber\\
&&p = t e^{-t}~sin(2\pi x)~sin(2\pi y).\nonumber
\end{eqnarray}
\end{example}
\noindent
In Table 2, we have shown the convergence
rates for backward Euler
method, respectively for $\bL^2$ and $\bH^1$-norms in velocity and $L^2$-norm
in pressure with $k=\mathcal{O}(h^{2})$. These results agree with the optimal theoretical
convergence rates obtained in
Theorem \ref{Tbd1}.
\small{
\begin{table}[ht!]
 \centering
 \begin{tabular}{|l|l|l|l|l|l|l|l|l|l|}
 \hline
  \small{$N$}&  h &\small{$\|\bu(t_n)-\bU^n\|$}  &\small{Rate}& \small{$\|\bu(t_n)-\bU^n\|_{\bH^1}$} & \small{Rate} & \small{$\|p(t_n)-P^n\|$}&\small{ Rate}  \\
             &    &                                       &                 &            &        &                  &   \\
 \hline
 \hline      & 1/4 &0.132916  &           & 3.736491 &          &0.989116&   \\
 \hline 577  & 1/8 &0.028166  & 2.238442  & 1.537594 & 1.281009 &0.186591& 2.406260  \\
 \hline 2433 & 1/16&0.003717  & 2.921475  & 0.463199 & 1.730969 &0.063322& 1.559099  \\
 \hline 9986 & 1/32&0.000473  & 2.971736  & 0.124017 & 1.901084 &0.016057& 1.979463 \\
 \hline 40449& 1/64&0.000063  & 2.894587  & 0.032022 & 1.953411 &0.006437& 1.318638 \\
 \hline
 \end{tabular}
 \vspace{.1cm}
 \caption{ Errors and convergence rates for backward Euler method with $k=\mathcal{O}(h^2)$. }
 \end{table}
}

\end{document}